\documentclass[10pt, USenglish, reqno]{article}

	\usepackage{amsmath}
	\usepackage{amssymb}
	\usepackage{amsthm}
	\usepackage[noend]{algpseudocode}
	\usepackage[USenglish]{babel}
	\usepackage{bm}
	\usepackage{comment}
	\usepackage{enumitem}
	\usepackage{etoolbox}
	\usepackage[T1]{fontenc}
	\usepackage[bottom=1.2in,left=1.1in,right=1.1in]{geometry}
	\usepackage{graphicx}
	\usepackage{lmodern}
	\usepackage{mathtools}
	\usepackage{nicefrac}
	\usepackage{subcaption}
	\usepackage[dvipsnames,table]{xcolor}
	\usepackage{xargs}
	\usepackage[%
		bookmarksdepth=4,
		breaklinks,
		hypertexnames = true,
		plainpages    = false,
		colorlinks = true,
		urlcolor   = Blue,
		linkcolor  = MidnightBlue,
		citecolor  = Green,
	]{hyperref}
	\usepackage{algorithm}
	\usepackage[capitalize,nameinlink]{cleveref}

	\graphicspath{{Pics/}}
	\theoremstyle{plain}
	\newtheorem{theorem}{Theorem}[section]
	\newtheorem{definition}[theorem]{Definition}
	\newtheorem{assumption}{Assumption}
		\Crefname{assumption}{Assumption}{Assumptions}
	\newtheorem{corollary}[theorem]{Corollary}
	\newtheorem{fact}[theorem]{Fact}
	\newtheorem{lemma}[theorem]{Lemma}
	\newtheorem{proposition}[theorem]{Proposition}

	\theoremstyle{definition}
	\newtheorem{remark}[theorem]{Remark}

	\renewcommand\theequation{\thesection.\arabic{equation}}
	\counterwithin{equation}{section}

	\setlist[enumerate]{leftmargin=*,widest*=3}
	\Crefname{figure}{Figure}{Figures}

	\newcommand{\numberthis}{\stepcounter{equation}\tag{\theequation}}

	\newlist{assenumerate}{enumerate}{1}
		\crefalias{assenumeratei}{assumption}
	\newlist{corenumerate}{enumerate}{1}
		\crefalias{corenumeratei}{corollary}
	\newlist{defenumerate}{enumerate}{1}
		\crefalias{defenumeratei}{definition}
	\newlist{factenumerate}{enumerate}{1}
		\crefalias{factenumeratei}{fact}
	\newlist{lemenumerate}{enumerate}{1}
		\crefalias{lemenumeratei}{lemma}
	\newlist{propenumerate}{enumerate}{1}
		\crefalias{propenumeratei}{proposition}
	\newlist{remenumerate}{enumerate}{1}
		\crefalias{remenumeratei}{remark}
	\newlist{thmenumerate}{enumerate}{1}
		\crefalias{thmenumeratei}{theorem}

	\newcommand{\coloneqq}{\mathrel{:=}}

	\DeclareMathOperator{\argmin}{argmin}
	
	\DeclareMathOperator{\DB}{D}
	\DeclareMathOperator{\dist}{dist}
	\DeclareMathOperator{\dom}{dom}
	\DeclareMathOperator{\epi}{epi}
	\DeclareMathOperator{\fix}{fix}
	\DeclareMathOperator{\graph}{gph}
	\DeclareMathOperator{\id}{id}
	\DeclareMathOperator{\jac}{J}
	\DeclareMathOperator{\lev}{lev}
	\DeclareMathOperator{\minimize}{minimize}
	\DeclareMathOperator{\prox}{prox}
	\DeclareMathOperator{\stt}{subject\ to}
	\DeclareMathOperator{\Ts}{t}		

	\newcommand{\C}{\mathcal{C}}
	\newcommand{\I}{\mathrm{I}}
	\newcommand{\N}{\mathbb{N}}
	\newcommand{\R}{\mathbb{R}}
	\newcommand{\Rinf}{\overline{\R}}
	\newcommand{\U}{\mathcal{U}}

	\newcommand{\FBE}{\Phi^{\hat{H}}}					
	\newcommand{\FBEG}{\phi^{H}}						
	\newcommand{\proxbreg}{\prox_{\phi}^h}
	\newcommand{\fbe}{\varphi^{\hat{h}}}				
	\newcommand{\Tt}{\operatorname{t}_{\hat{h}}}		
	\newcommand{\Res}{\operatorname{R}_{\hat{H}}}		
	\newcommand{\res}{\operatorname{r}_{\hat{h}}}
	\newcommand{\MM}{\mathcal{M}_{\hat{H}}}

	\newcommand{\tbz}[2][1=i]{\tilde{\bm z}^{#2}_{#1}}

	\newcommandx{\tz}[3][1=i,3={}]{\bar{\bm z}^{#2}_{#1}}
	\newcommandx{\ttz}[3][1=i,3={}]{\tilde{z}^{#2}_{#1}}

	\newcommand{\vect}[1]{\boldsymbol{#1}}

	\makeatletter
		\renewcommand{\alglinenumber}[1]{{%
			\footnotesize
			\textbf{\thealgorithm}.%
			\rlap{\oldstylenums{\arabic{ALG@line}}:}%
			\hphantom{\oldstylenums{18}:}%
		}}
		\renewcommand{\theALG@line}{\thealgorithm.\oldstylenums{\arabic{ALG@line}}}

		\newcounter{algorithmicH}
		\let\OLDalgorithmic\algorithmic
		\renewcommand{\algorithmic}{%
			\setcounter{ALG@line}{0}%
			\stepcounter{algorithmicH}%
			\OLDalgorithmic%
		}%
		\newcommand{\theHALG@line}{ALG@line.\thealgorithmicH.\arabic{ALG@line}}%
	\makeatother

\begin{document}

	\title{%
		SPIRAL: A superlinearly convergent incremental proximal algorithm for nonconvex finite sum minimization%
		\thanks{%
			P. Behmandpoor and M. Moonen acknowledge the research work carried out at the ESAT Laboratory of KU Leuven, in the frame of Research Project FWO nr. G0C0623N 'User-centric distributed signal processing algorithms for next generation cell-free massive MIMO based wireless communication networks' and Fonds de la Recherche Scientifique - FNRS and Fonds voor Wetenschappelijk Onderzoek - Vlaanderen EOS Project no 30452698 '(MUSE-WINET) MUlti-SErvice WIreless NETworks'.
			The scientific responsibility is assumed by its authors.
			The work of P. Latafat was supported by the
			Research Foundation Flanders (FWO) grants 1196820N and 12Y7622N.
			The work of P. Patrinos was supported by the
			Research Foundation Flanders (FWO) research projects G0A0920N, G086518N, G086318N, and G081222N;
			Research Council KU Leuven C1 project No. C14/18/068;
			Fonds de la Recherche Scientifique -- FNRS and the Fonds Wetenschappelijk Onderzoek -- Vlaanderen under EOS project 30468160 (SeLMA);
			European Union's Horizon 2020 research and innovation programme under the Marie Sk\l odowska-Curie grant agreement No. 953348.%
			The work of A. Themelis was supported by the Japan Society for the Promotion of Science (JSPS) KAKENHI grant JP21K17710.%
		}%
	}

	\author{%
		Pourya Behmandpoor\thanks{%
			KU Leuven, Department of Electrical Engineering (ESAT-STADIUS), Kasteelpark Arenberg 10, 3001 Leuven, Belgium.
			\textit{Emails:}
			\textsf{%
				\{%
					\href{mailto:pourya.behmandpoor@esat.kuleuven.be}{pourya.behmandpoor},%
					\href{mailto:puya.latafat@esat.kuleuven.be}{puya.latafat},%
					\href{mailto:marc.moonen@esat.kuleuven.be}{marc.moonen},%
					\href{mailto:panos.patrinos@esat.kuleuven.be}{panos.patrinos}%
				\}%
				\href{mailto:pourya.behmandpoor@esat.kuleuven.be,puya.latafat@esat.kuleuven.be,marc.moonen@esat.kuleuven.be,panos.patrinos@esat.kuleuven.be}{@esat.kuleuven.be}%
			}%
		}%
		\and Puya Latafat\footnotemark[2]%
		\and Andreas Themelis\thanks{%
			Kyushu University, Faculty of Information Science and Electrical Engineering (ISEE), 744 Motooka Nishi-ku, 819-0395 Fukuoka, Japan.
			\textit{Email:}
			\textsf{%
				\href{mailto:andreas.themelis@ees.kyushu-u.ac.jp}{andreas.themelis@ees.kyushu-u.ac.jp}%
			}%
		}%
		\and Marc Moonen\footnotemark[2]%
		\and Panagiotis Patrinos\footnotemark[2]%
	}

	\date{}

	\maketitle

	\begin{abstract}
		We introduce SPIRAL, a {\bf S}u{\bf P}erlinearly convergent {\bf I}ncremental p{\bf R}ox\-i\-mal {\bf AL}\-go\-rithm, for solving nonconvex regularized finite sum problems under a relative smoothness assumption.
		Each iteration of SPIRAL consists of an inner and an outer loop.
		It combines incremental gradient updates with a linesearch that has the remarkable property of never being triggered asymptotically, leading to superlinear convergence under mild assumptions at the limit point.
		Simulation results with L-BFGS directions on different convex, nonconvex, and non-Lipschitz differentiable problems show that our algorithm, as well as its adaptive variant, are competitive to the state of the art.
	\end{abstract}





	\section{Introduction}

		We study nonconvex nonsmooth finite sum optimization problems of the form:
		\begin{equation}\label{eq:problem_formulation}
			\textstyle \minimize_{z \in \R^n} \varphi(z)  \coloneqq f(z)+g(z),\quad
			\text{where} \; f(z) \coloneqq\tfrac{1}{N}\textstyle\sum_{i=1}^{N} f_i(z).
		\end{equation}
		The following basic assumptions are considered throughout the paper:

		\begin{assumption}[basic assumptions]\label{ass:basic}~
		    \begin{assenumerate}[label=A\arabic*,ref=\ref*{ass:basic}.A\arabic*,itemsep=0pt,topsep=0pt]
			\item\label{ass:basic:1}%
				$f_i: \R^n \to \R $ is $L_i$-smooth relative to a distance-gen\-er\-at\-ing function $h_i:\R^n\to\R$ (cf. \cref{def:relative_smoothness,def:legendre_kernel}), $i \in [N]\coloneqq\{1,\dots,N\}$;
			\item
				$g: \R^n \to \Rinf $ is proper and lower semicontinuous (lsc);
			\item a solution exists: $\arg\min \varphi \neq \emptyset$.
		    \end{assenumerate}
		\end{assumption}

		The minimization problem \eqref{eq:problem_formulation} has gained considerable attention across various disciplines including machine learning (ML), signal and image processing, statistics, and control.
		Despite an upsurge in developing optimization methods to address such a problem, the potential of low-memory quasi-Newton methods has largely been neglected which can be partially attributed to the absence of theoretical foundations for handling nonsmooth settings.
		In the smooth strongly convex settings, competitive convergence rates compared to typical ML approaches have been documented in the ML community \cite{moritz2016linearly}. This work aims to address such large-scale problems in their full generality in the nonconvex, nonsmooth problem settings.

		Stochastic gradient descent (SGD) is commonly employed for finite sum minimization problems.
		Despite it involving simple iterations, SGD requires a diminishing stepsize and, even in the strongly convex setting, can only achieve sublinear rates of convergence. These limitations have prompted the development of several stochastic and incremental methods such as SAG \cite{schmidt2017minimizing}, SAGA \cite{defazio2014saga}, SDCA \cite{shalevshwartz2013stochastic}, SPIDER \cite{fang2018spider}, SVRG \cite{johnson2013accelerating} and its extensions \cite{reddi2016stochastic,ge2019stabilized}, SARAH \cite{nguyen2017sarah}, and zeroSARAH \cite{li2021zerosarah}, which primarily target smooth functions ($g=0$) and are often restricted to the convex regime.
		To accommodate composite nonsmooth cost functions of the form \eqref{eq:problem_formulation}, studies such as \cite{cai2022cyclic}, proxSGD \cite{ghadimi2016mini}, proxSAGA and proxSVRG \cite{reddi2016proximal}, proxSARAH \cite{pham2020proxsarah}, and SpiderBoost \cite{wang2018spiderboost} have emerged recently.

		The majority of the methods mentioned above incorporate a combination of stochastic and deterministic components in addressing the finite sum problem, aiming to diminish the variance of iterates toward the optimal point.
		Notably, algorithms such as SAGA, SVRG, and SARAH employ an outer loop to incorporate full gradients as the deterministic enhancement, along with an inner loop that incorporates stochastic gradients using randomized sampling with replacement.
		Furthermore, these algorithms adopt fixed stepsizes, in contrast to SGD which necessitates diminishing stepsizes to mitigate variance.
		In line with the spirit of these methods, the porposed algorithm also utilizes both inner and outer loops.

		In its inner loop, the algorithm investigated in this study can be perceived as an incremental approach with a (shuffled) cyclic (randomized without replacement) sweeping rule.
		It should be mentioned that when combined with SGD, this sweeping rule demonstrates superior convergence and implementation efficiency \cite{recht2013parallel,bengio2012practical} compared to the randomized sweeping rule with replacement.
		Moreover, the analysis of SGD with (randomized) sampling without replacement has recently emerged in convex regimes, providing enhanced bounds compared to standard SGD \cite{cai2023empirical,cha2023tighter,mishchenko2020random,gurbuzbalaban2021random,haochen2019random}.
		Beyond SGD, the proposed algorithm, in its basic form without a linesearch, can be considered as a memory-efficient variant of Finito \cite{defazio2014finito} and MISO \cite{mairal2015incremental}.
		It is worth noting that DIAG \cite{mokhtari2018surpassing}, proposed independently, studies the Finito/MISO algorithm under a cyclic sweeping rule and in the strongly convex case.
		More recently, \cite{latafat2021block} provided a comprehensive study of the above algorithms in the fully nonconvex setting.
		However, the above are all limited to first-order methods.

		In its outer loop, one distinguishing characteristic of the proposed algorithm, which sets it apart from stochastic algorithms like SVRG and SARAH, is its utilization of quasi-Newton directions integrated with a linesearch while preserving the advantageous low-memory characteristic.
		In the context of this study, various methodologies have been explored to address nonconvex nonsmooth composite functions by employing quasi-Newton directions. For instance, methods presented in \cite{themelis2018forward,themelis2019supermann,ahookhosh2021bregman} have demonstrated the application of quasi-Newton directions to achieve superlinear convergence rates, albeit limited to scenarios involving a single smooth function within the composite cost.
		In the finite sum setting, approaches proposed by \cite{moritz2016linearly} and \cite{sadeghi2021hybrid} have utilized quasi-Newton updates with global convergence guarantees and linear convergence rates.
		Furthermore, \cite{zhang2022variance} has extended the utilization of quasi-Newton directions to decentralized learning scenarios.

		To attain a superlinear convergence rate, the IQN method \cite{mokhtari2018iqn} has integrated quasi-Newton directions with incremental updates, albeit with only local convergence guarantees. Conversely, the approach introduced in \cite{rodomanov2016superlinearly} also exhibits a superlinear convergence rate but necessitates Hessian evaluation. It is noteworthy that the aforementioned algorithms are applicable in (strongly) convex cases. However, within the nonconvex nonsmooth setting, the algorithm proposed by \cite{yang2021stochastic} stands out with global convergence guarantees when the nonsmooth term is convex.

		One of the restrictive aspects of the aforementioned works is that the cost functions are Lipschitz differentiable.
		However, in numerous practical applications, although the cost functions are differentiable, they fail to satisfy Lipschitz continuity assumptions on their gradients.
		This issue is exemplified in \cref{sec:phaseRetrieval}.
		To tackle such cost functions, the proposed algorithm goes beyond the classical notion of smoothness and employs the concept of relative smoothness, as introduced in \cite{bauschke2017descent,lu2018relatively}.
		In connection with the proposed method, stochastic mirror descent (SMD) methods incorporate relative smoothness. Notable references in this area include \cite{beck2003mirror,nedic2014stochastic,hanzely2021fastest,davis2018stochastic}.
		Within the convex setting, PLIAG \cite{zhang2021proximal} has been introduced as the Bregman variant of IAG \cite{bertsekas2000gradient,blatt2007convergent,vanli2018global}. Furthermore, \cite{dragomir2021fast} explores Bregman stochastic gradient descent (BSGD). In the nonconvex regime, \cite{latafat2021bregman} investigates a Bregman variant of Finito/MISO.

		While the literature often assumes convexity for the nonsmooth term $g$, our proposed method, as indicated in \cref{ass:basic}, allows the nonconvex nature of this term. This enables the algorithm to effectively handle a wide range of nonconvex constraints, including rank constraints and $\ell_0$-norm ball constraints, as well as nonconvex regularizers such as $\ell^p$ with $p \in [0,1)$.

		Motivated by the aforementioned advancements and recognizing the existing limitations in the literature, the proposed method addresses the optimization of regularized nonsmooth nonconvex cost functions, allowing the gradients of differentiable functions in the finite sum to be non-Lipschitz.
		To the best of our knowledge, none of the currently available methods in the literature that exhibit superlinear convergence rates have explicitly addressed the challenge of handling non-Lipschitz differentiable functions within nonsmooth nonconvex finite sum settings.

		\subsection*{Contributions}

			The main contributions of the paper are as follows:
			\begin{enumerate}[%
				leftmargin=*,
				label={\arabic*.},
				widest*=3,
			]
			\item
				We propose SPIRAL with convergence guarantees for a wide class of finite sum problems. Not only are both the nonsmooth regularizer \(g\) and the finite sum terms \(f_i\) all allowed to be nonconvex, but also \(f_i\) functions do not need to have Lipschitz-continuous gradients. Moreover, unlike Finito/MISO/DIAG, SPIRAL requires only \(\mathcal{O}(n)\) memory allocation.
			\item
				When the nonsmooth term is convex, we show that SPIRAL converges superlinearly when the employed quasi-Newton directions are superlinear (cf. \cref{def:superlinearDirection}) and
				the linesearch will eventually never be invoked (cf. \cref{thm:superlinear}) under mild assumptions.
				This is also supported by our simulation results.
				Moreover, global (as opposed to local) convergence is guaranteed regardless of any assumptions placed on the quasi-Newton directions or the convexity of the nonsmooth term (cf. \cref{lem:descent_lemma}).
			\item
				Finally, an adaptive variant employing appropriate backtracking linesearch is introduced that adapts to the local relative smoothness moduli of \(f_i\) while maintaining convergence guarantees and the \(\mathcal{O}(n)\) memory requirement.
			\end{enumerate}

	\section{Preliminaries}\label{sec:setup}

		\subsection{Notation}\label{sec:auxiliary:not}

			In this section, we provide basic notations. The interested reader may refer to \cite{rockafellar1970convex,rockafellar2009variational} for details.
			The set of natural numbers is denoted by \(\N = \{0,1,2,\ldots\}\).
			The set of real and extended-real numbers are \(\R\coloneqq(-\infty,\infty)\) and \(\Rinf\coloneqq\R\cup\{\infty\}\), and the set of positive reals is denoted by \(\R_+\coloneqq [0, \infty)\).
			We also use the notation $[N] \coloneqq \{1,2,\cdots,N\}$.
			We denote by \(\langle{}\cdot{},{}\cdot{}\rangle\) and \(\|{}\cdot{}\|\) the standard Euclidean inner product and the induced norm.
			The distance of a point $x\in\R^n$ to a nonempty set $S\subseteq\R^n$ is given by $\dist(x,S)=\inf_{z\in S}{\|z-x\|}$.
			For a vector $\bm w=(w_1,\ldots,w_r)\in\R^{\sum_i n_i}$, $w_i\in\R^{n_i}$ is used to denote its $i$-th block coordinate.
			The identity operator is denoted by \(\id\).

			For a sequence \((x^k)_{k\in\N}\) we write \((x^k)_{k\in\N}\subseteq E\) to indicate that \(x^k\in E\) for all \(k\in\N\). We use the following notions of convergence rate: a sequence $(x^k)_{k\in\N}$ is said to converge to a point $x^\star$:
			\begin{itemize}
				\item (at least) \emph{$Q$-linearly} (with quotient rate)  with $Q$-factor given by $\sigma\in(0,1)$, if there exists $k_0 \in \N$ such that for all $k \geq k_0$,
				\[
					\|x^{k+1}-x^\star\|\leq\sigma\|x^k-x^{\star}\|.
				\]

				\item (at least) \emph{R-linearly} (with root rate) if there exists a sequence of nonnegative scalars $({v^k})_{k\in\N}$ such that  $\|x^k-x^\star\|\leq v^k$ and $(v^k)_{k\in\N}$ converges Q-linearly to zero.

				\item superlinearly if either $x^k=x^\star$ for some $k\in \N$ or
				\[
					\lim_{k\to \infty}\frac{\|x^{k+1}-x^\star\|}{\|x^k-x^\star\|} = 0.
				\]
			\end{itemize}


			We use the notation \(Q:\R^n\rightrightarrows\R^m\) to indicate a mapping from each point $x\in\R^n$ to a subset $Q(x)$ of $\R^m$.
			The \emph{graph} of \(Q\) is the set
			\(
				\graph Q
			\coloneqq
				\{(x,y)\in\R^n\times\R^m:y\in Q(x)\}
			\),
			and the set of its \emph{fixed points} is defined as
			\(
				\fix Q
			\coloneqq
				\{x\in\R^n:x\in Q(x)\}
			\).
			We say that \(Q\) is \emph{outer semicontinuous (osc)} if \(\graph Q\) is a closed subset of \(\R^n\times\R^m\), and \emph{locally bounded} if for every bounded \(U\subset\R^n\) the set \(\bigcup_{x\in U}Q(x)\) is bounded.

			The \emph{domain} of an extended-real-valued function \(\phi:\R^n\to\Rinf\) is the set
			\(
				\dom \phi
			\coloneqq
				\{x\in\R^n:\phi(x)<\infty\}
			\)
			and
			\(
				\epi \phi
			\coloneqq
				\{(x,\alpha)\in\R^n\times\R:\phi(x)\leq\alpha\}
			\)
			is its \emph{epigraph} set.
			Function \(\phi\) is said to be \emph{proper} if \(\dom \phi\neq\emptyset\), and \emph{lower semicontinuous (lsc)} if \(\epi \phi\) is a closed subset of \(\R^{n+1}\).
			We say that \(\phi\) is \emph{level bounded} if its \(\alpha\)-sublevel set
			\(
				\lev_{\leq\alpha}\phi
			\coloneqq
				\{x\in\R^n:\phi(x)\leq\alpha\}
			\)
			is bounded for all \(\alpha\in\R\).
			The indicator function $\delta_X$  of a nonempty set $X\subseteq \R^n$ is defined by
			\begin{equation}
				\delta_X(x) \coloneqq \begin{cases}
				0  & \ \ \ \  \textrm{if} \;\, x\in X\\
				+\infty  & \ \ \ \ \textrm{if} \;\, x\notin X.
			\end{cases}
			\end{equation}
			The indicator function $\delta_X$ is closed if and only if $X$ is a closed set.

			We denote by \(\hat\partial\phi:\R^n\rightrightarrows\R^n\) the \emph{regular sub-differential} of \(\phi\), where
			\begin{equation*}
				v\in\hat\partial \phi(\bar x)
			\iff
				\liminf_{\bar x\neq x\to\bar x}{
					\frac{\phi(x)-\phi(\bar x)-\langle v,x-\bar x\rangle}{\|x-\bar x\|}
				}
			\geq
				0.
			\end{equation*}
			The regular sub-differential is closed- and convex-valued.
			The (limiting) \emph{sub-differential} of \(\phi\) is \(\partial\phi:\R^n\rightrightarrows\R^n\), where
			\(
				v\in\partial \phi(x)
			\)
			iff \(x\in\dom \phi\) and there exists a sequence \((x^k,v^k)_{k\in\N}\subseteq\graph\hat\partial \phi\) such that
			\(
				(x^k,\phi(x^k),v^k)
			\to
				(x,\phi(x),v)
			\)
			as \(k\to\infty\).

			A necessary condition for local minimality of \(x\) for \(\phi\) is \(0\in\hat\partial \phi(x)\), see \cite[Thm. 10.1]{rockafellar2009variational}.
			Finally, the set of $r$ times continuously differentiable functions over $\R^n$ is denoted by $\C^r =\C^r(\R^n)$.

		\subsection{Relative smoothness}

			We start by formally defining the notion of distance generating function, Bregman distance, and relative smoothness \cite{bauschke2017descent,lu2018relatively,teboulle2018simplified}.

			\begin{definition}[Distance-generating function (dgf)]\label{def:legendre_kernel}
				A strictly convex function $h: \R^n \to \R$ that is continuously differentiable everywhere will be referred to as a distance-generating function (dgf).
			\end{definition}

			While some results here can be presented with $\dom h \subset \R^n$, for the sake of simplicity and global convergence analysis, we continue with $\dom h = \R^n$. The Bregman distance associated with a dgf is defined as:

			\begin{definition}[Bregman distance]\label{def:bregman_distance}
				Given a dgf $h: \R^n \to \R$, the Bregman distance $\DB_h: \R^n \times \R^n \to \R_+$ is defined as,
				\begin{equation}
					\DB_h(y,x) \coloneqq h(y) - h(x) - \langle \nabla h(x), y-x\rangle \quad \text{for all } x,y \in \R^n.
					\label{eq:bregman}
				\end{equation}
			\end{definition}

			\begin{definition}[Relative smoothness \protect{\cite[Def. 2.2]{bolte2018first}}]\label{def:relative_smoothness}
				A function $f: \R^n \to \R$ is smooth relative to a dgf $h:\R^n \to \R$ if there exists $L \geq 0$ such that $Lh \pm f$ are convex functions on $\R^n$.
				In this case, we say that \(f\) is \(L\)-smooth (relative to \(h\)) to make the modulus \(L\) explict.
			\end{definition}

			\begin{fact}[Descent lemma \protect{\cite[Lem. 2.1]{bolte2018first}}]\label{fact:smoothness}
				If $f:\R^n \to \R$ is $L$-smooth relative to a dgf $h:\R^n \to \R$, then for all \(x,y\in\R^n\)
				\[ 
					|f(y) - f(x) - \langle \nabla f(x), y-x\rangle| \leq L \DB_h(y,x).
				\]
			\end{fact}
			Note that in the Euclidean case, the dgf and the corresponding Bregman distance reduce to $h=\tfrac{1}{2}\|\cdot\|^2$ and $\DB_h(y,x)=\tfrac{1}{2}\|y-x\|^2$, respectively, and \cref{fact:smoothness} reduces to the ordinary descent lemma \cite[Prop. A.24]{bertsekas2016nonlinear} for smooth functions.

			Relative to a dgf $h:\R^n \to \R$, the \emph{(left) Bregman proximal mapping} of a proper lsc function $\phi: \R^n \to \Rinf$ is the set-valued mapping $\proxbreg:\R^n\rightrightarrows\R^n$ defined as \cite[Def. 2.2]{kan2012moreau}
			\begin{align}
				\label{eq:prox_def}
				\proxbreg(x)
			\coloneqq{} &
				\argmin_{w \in \R^n}\bigl\{\phi(w) + \DB_h(w,x)\bigr\},
			\shortintertext{and its value function is the \emph{Bregman Moreau envelope} $\phi^h: \R^{n} \to \R$}
				\label{eq:prox_def_val}
				\phi^h(x)
			\coloneqq{} &
				\min_{w \in \R^n}\bigl\{\phi(w) + \DB_h(w,x)\bigr\}.
			\end{align}
			Moreover it is evident from \eqref{eq:prox_def} and \cref{def:bregman_distance} that if $v \in \proxbreg(x)$, then
			\begin{equation}\label{thm:OC}
				\nabla h(x) - \nabla h(v) \in \hat\partial \phi(v),
			\end{equation}
			and the converse also holds when $\phi$ is convex.
			Whenever the superscript $h$ is omitted from $\proxbreg$, it refers to the \emph{Euclidean proximal mapping} with $h=\tfrac{1}{2}\|\cdot\|^2$.

	\section{Proposed algorithm}\label{sec:Alg}

		\begin{algorithm}[t]
			\caption{SPIRAL}
			\label{alg:LBFinito}

				\begin{algorithmic}[1]
		\item[\textbf{Initialize}]
			\begin{tabular}[t]{@{}l@{}}
				\(z^{\rm init}\in\R^n\),~
				\(\gamma_i\in(0,\nicefrac{N}{L_{i}})\),~%
				\(\forall i\in[N]\),~
				\(s^{\rm 0} = \sum_{i=1}^N{\frac{1}{\gamma_i}\nabla h_i(z^{\rm init}) - \frac{1}{N} \nabla f_i(z^{\rm init})}\),~
			\\
				maximum number of backtracks \(q_{\rm max} \in \N\cup\{\infty\}\) {\scriptsize{(e.g. \(q_{\rm max}=2\))}},~
				\(\beta \in (0,1)\)
			\end{tabular}

		\item[\textbf{Repeat for} $k=0,1,\dots,K$]
			\State\label{alg:1:t1}
				\(z^k\in\Ts(s^k)\)

			\State\label{alg:1:s2}
			\(
				\bar s^k
				=
				\sum_{i=1}^N{\frac{1}{\gamma_i}\nabla h_i(z^k) - \frac{1}{N} \nabla f_i(z^k)}
			\)
			\Comment{(full update)}%

			\State\label{alg:1:t2}
				\(v^k\in\Ts(\bar s^k)\)

			\State \label{alg:1:dir}%
				choose $d^k\in\R^n$ at $z^k$ {\footnotesize(e.g. based on a quasi-Newton method for solving $\res(z)=0$)}

			\State \label{alg:1:ls}%
				set $\tau_k = 1, q_k=0$
			\Comment{(linesearch)}
				\begin{enumerate}[label=\alph*:,ref=\ref*{alg:1:ls}.\alph*,widest*=4]
					\item\label{alg:1:ls:1} \(u^k=\tau_k z^k+(1-\tau_k)v^k+\tau_k{d^k}\)

					\item\label{alg:1:ls:2}
						\(
							\tilde s^{k}
							=
							\sum_{i=1}^N{\frac{1}{\gamma_i}\nabla h_i(u^k) - \frac{1}{N} \nabla f_i(u^k)}
						\)
					\Comment{(full update)}

					\item \label{alg:1:ls:3}%
						\(y^k\in\Ts(\tilde s^{k})\)

					\item \label{eq:ls_condition}%
						\textbf{if} \(\mathcal{L}(y^k,u^k)\leq\mathcal{L}(v^k,z^k)\)

						\quad go to step \ref{alg:1:sts}%
					\item
						\textbf{else if} \(q_k = q_{\max}\) \textbf{then}

						\quad
						\(u^k = v^k\),~
						\(\tilde s^{k}=\sum_{i=1}^N{\frac{1}{\gamma_i}\nabla h_i(u^k) - \frac{1}{N} \nabla f_i(u^k)}\),
						and go to step \ref{alg:1:sts}

					\item
						\textbf{else}

						\quad
						\(\tau_k \gets \beta \tau_k\),~
						\(q_k \gets q_k + 1\),
						and go to step \ref{alg:1:ls:1}%
				\end{enumerate}
			\State \label{alg:1:sts}%
				\(s^k\gets\tilde s^k\)

			\For{\(\ell = 1,\ldots, N\)}\label{alg:forloop}%
			\Comment{(incremental loop)}

			\State \label{alg:1:rnd}%
				randomly choose \(i^\ell\in[N]\) without replacement

			\State \label{alg:1:t3}%
				\(\tilde z_{i^\ell}^k\in\Ts(s^k)\)

			\State \label{alg:forloop:end}%
				\(
					s^k
				\gets
					s^k +\left[
						\frac{1}{\gamma_{i^\ell}}\nabla h_{i^\ell}(\tilde z^k_{i^\ell}) - \frac{1}{N} \nabla f_{i^\ell}(\tilde z^k_{i^\ell})
						-
						\frac{1}{\gamma_{i^\ell}}\nabla h_{i^\ell}(u^k) + \frac{1}{N} \nabla f_{i^\ell}(u^k)
					\right]
				\)
			\EndFor

		\State
			$s^{k+1} \gets s^k$

		\item[\textbf{Return} $z^K$]
		\end{algorithmic}
		\end{algorithm}

		The proposed method, SPIRAL, is outlined in \cref{alg:LBFinito} to address the optimization problem \eqref{eq:problem_formulation}.
		SPIRAL employs the set-valued mapping $\Ts:\R^n\rightrightarrows\R^n$ as its major oracle, defined as
		\begin{equation}\label{eq:t}
			\textstyle \Ts(s)
		\coloneqq
			\argmin_{w\in\R^n}\left\{
				g(w)
				+
				\sum_{i=1}^N\tfrac{1}{\gamma_i}h_i(w)
				-
				\langle s,w\rangle
			\right\},
		\end{equation}
		where $h_i$ is the dgf corresponding to function $f_i$ as in \cref{ass:basic:1}.
		Within both the outer and inner loops, the mapping $\Ts$, utilized in steps \ref{alg:1:t1}, \ref{alg:1:t2}, \ref{alg:1:ls:3}, and \ref{alg:1:t3}, represents the proximal steps.
		It is important to note that in the case of Euclidean space, where the functions $f_i$ are $L_i$-smooth relative to $h_i = \frac{1}{2}\|\cdot\|^2$, the oracle $\Ts$ reduces to the Euclidean proximal mapping with updates of the form $z^k \in \prox_{\hat{\gamma} g}(\hat{\gamma}s^k)$, where $\hat{\gamma}^{-1} = \sum_{i=1}^N \gamma_i^{-1}$.
		Similarly, the iterates $v^k$, $y^k$, and $\ttz[i]{k}$ are updated using the same $\prox_{\hat{\gamma} g}$ function.
		A detailed description of the Euclidean version of the algorithm can be found in Appendix \ref{sec:Euclidean}.

		SPIRAL employs the function $\mathcal{L}: \R^n \times \R^n \to \R$ defined by
		\begin{equation}\label{eq:Lyapunov}
			\mathcal{L}(y,x) \coloneqq \varphi(y)+\DB_{\hat{h}}(y,x) \quad \text{with} \;\, y \in \Ts(\nabla \hat{h}(x))
		\end{equation}
		in its linesearch,
		where $\hat{h}_i: \R^n \to \R$ and $\hat{h}:\R^n \to \R$ are
		\begin{equation}\label{eq:dgf}
			\textstyle \hat{h}_i \coloneqq \frac{1}{\gamma_i}h_i - \frac{1}{N} f_i, \quad \hat{h} \coloneqq \sum_{i=1}^N \hat{h}_i.
		\end{equation}
		The function $\mathcal{L}$ is considered as a suitable Lyapunov function in our convergence analysis.
		The linesearch in step \ref{alg:1:ls} interpolates the iterate $u^k$ between the candidate fast update $z^k+d^k$, corresponding to \(\tau_k=1\), and the safeguard step $v^k$ which is approached as \(\tau_k\searrow0\). The iterate $u^k$ is selected whenever it is a descent direction for the Lyapunov function $\mathcal{L}$ (cf. \cref{rem:wellDefinedness} regarding the well definedness of the linesearch).

		One distinguishing characteristic of SPIRAL, which sets it apart from stochastic algorithms such as SVRG and SARAH, is its utilization of directions $d^k$ in step \ref{alg:1:dir} based on second-order-like information of the set-valued residual mapping $\res:\R^n\rightrightarrows\R^n$ defined as
		\begin{equation}\label{eq:FB_residual}
			\res \coloneqq \id - \Ts \circ \nabla \hat{h}.
		\end{equation}
		This feature allows SPIRAL to achieve a superlinear convergence rate, given certain mild conditions (cf. \cref{thm:superlinear}).
		Given that the inclusion $0\in \res(z^\star)$ holds at a stationary point $z^\star$ (as discussed in \cref{sec:reform} and expressed in \eqref{eq:fixed_point}), the direction $d^k$ is determined based on solving this inclusion. Semismooth Newton directions \cite{themelis2019supermann} can be utilized to compute $d^k$; however, this approach relies on access to second-order oracle information. As an alternative to circumvent this requirement, quasi-Newton methods can be employed to compute the directions $d^k$.
		Specifically, we employ the update rule
		\begin{equation}
			d^k = -H^k \res(z^k)
			\label{eq:update_rule_of_d}
		\end{equation}
		where $H^k$ represents a linear operator that approximates the second-order information of the residual mapping $\res$.
		It is worth noting that although $\res$ is a set-valued mapping, it typically exhibits single-valuedness and other desirable properties in the vicinity of stationary points of the objective function $\varphi$ (see e.g. \cite[Thm. 3.10 and 3.11]{ahookhosh2021bregman}).
		Consequently, the updates for $H^k$ can be implemented using popular quasi-Newton methods such as Broyden's method, BFGS, and L-BFGS (cf. \cref{sec:spl} for further discussion).

		SPIRAL can accommodate various sweeping rules depending on the memory requirements. The following remark comments on two different settings.

		\begin{remark}[sweeping rule in the incremental loop]\label{rem:incremental}~
			\begin{remenumerate}[label=(\roman*),ref=\ref*{rem:incremental}(\roman*),widest*=2]
			\item\emph{low-memory setting:}
				in this setting, SPIRAL employs a (shuffled) cyclic (randomized without replacement) sweeping rule within the incremental loop, and unlike Finito/MISO/DIAG methods, it does not require storing individual function gradients $\nabla f_i$. Instead, SPIRAL only requires storing the finite sum of gradients as $\nabla \hat{h}$, using the vectors $s^k$, $\tilde s^k$, and $\bar s^k$.
				Additionally, in the incremental loop, updating this sum in the vector $s^k$ is carried out in a memory-efficient manner since the vector $u^k$ is known and fixed due to choosing $i^\ell$ without replacement. This advantageous characteristic allows SPIRAL to require only \(\mathcal{O}(n)\) memory allocation, making it suitable for large-scale optimization problems.
			\item\emph{high-memory setting:}
				in this setting, by utilizing a memory of $\mathcal{O}(nN)$ and saving $\tilde{z}^k_j$ for all $j\in[N]$, like Finito/MISO/DIAG, the incremental loop can be replaced by a randomized loop of an arbitrary depth.
			\end{remenumerate}
		\end{remark}

		\noindent\textbf{Computational complexity.}
			The overall computational complexity, measured in terms of gradient evaluations per iteration, is $\mathcal{O}\big((4+\kappa)N\big)$. This includes
			$2+\kappa$ full gradient evaluations performed outside of the inner loop, where $\kappa$ represents the number of backtracks in the linesearch, and two gradient evaluations performed in each iteration of the inner loop.
			It is worth noting that in certain problems like least squares, nonnegative principal component analysis, and logistic regression, the gradients can be obtained by storing the inner product between data points and the evaluated points. Consequently, the gradient evaluations $\nabla \hat{h}_{i}(u^k)$ in step \ref{alg:forloop:end} can be derived from the computations in step \ref{alg:1:ls:2}, thereby reducing the computational complexity to $\mathcal{O}\big((3+\kappa)N\big)$.
			In addition, based on numerical experiments, it is beneficial to limit the number of backtracks using a maximum value $q_{\max}$.
			It is important to note that, as demonstrated in \cref{thm:superlinear}, under mild conditions, the unit stepsize $\tau_k = 1$ is eventually always accepted ($\kappa=0$), allowing for pure quasi-Newton type updates in step \ref{alg:1:ls:1} and avoiding any further backtracks and computation.

	\section{Convergence Analysis}\label{sec:conv}
		In the subsequent subsections, we study the convergence of SPIRAL by reformulating the problem \eqref{eq:problem_formulation} in the lifted space.
		We will thus recast the Lyapunov function $\mathcal{L}$ in \eqref{eq:Lyapunov} into this new space and study it along the iterates generated by SPIRAL.
		In the lifted space, we can show that block coordinate updates in the incremental loop result in sufficient descent for the Lyapunov function, while this is not necessarily the case for the objective function.
		Utilizing newly reformulated operators, we then outline \cref{alg:LBFinito} in the lifted space.
		Having new insight into the mechanism of SPIRAL, in \crefrange{sec:subseq}{sec:glb} we establish its convergence in various regimes.

		\subsection{Problem Reformulation}\label{sec:reform}
			We recast \eqref{eq:problem_formulation} as the following lifted consensus optimization problem:
			\begin{align*}
				\minimize_{\vect z=(z_1,\dots,z_N) \in \R^{Nn}} \Phi(\vect z)
			\coloneqq{} &
				\underbracket[0.5pt]{\frac{1}{N}\sum_{i=1}^{N} f_i(z_i)}_{\coloneqq F(\bm z)}
				+
				\underbracket[0.5pt]{\frac{1}{N}\sum_{i=1}^{N} g(z_i) + \delta_\Delta(\vect z)}_{\coloneqq G(\bm z)}
			\stepcounter{equation}\tag{\theequation}
			\label{eq:recasted_problem_formulation}
			\\
			\text{where}\quad
				\Delta
			\coloneqq{} &
				\left\{\vect z = (z_1,\dots,z_N)\in \R^{Nn}|z_1=\dots=z_N\right\}
			\end{align*}
			is the \emph{consensus set}.
			Note that $\Phi(\vect z) = \varphi(z)$ whenever $\vect z=(z,\dots,z) \in \Delta$.
			Define
			\begin{equation}\label{eq:H_hat_def}
				\MM(\vect w, \vect z) \coloneqq \Phi(\vect w) + \DB_{\hat{H}}(\vect w, \vect z), \quad
				\text{with  } \; \hat{H}(\vect x) \coloneqq \sum_{i=1}^N \hat{h}_i(x_i),
			\end{equation}
			for any $(\vect w, \vect z) \in \R^{Nn} \times \R^{Nn}$, where $\hat{h}_i$ is defined in \eqref{eq:dgf}.
			The model $\MM$ is in particular a \emph{majorizing model} of $\Phi$, in that from \eqref{eq:dgf}, and \cref{def:legendre_kernel,def:relative_smoothness}, it is evident that whenever $\gamma_i\in (0, \nicefrac{N}L_i)$, $\hat{h}_i$ and consequently $\hat{H}$ are dgfs, hence $\DB_{\hat{H}}(\vect w, \vect z) \geq 0$ for all $\vect w, \vect z \in \R^{Nn}$ and%
			\begin{enumerate}[label=(\roman*),widest=2]
			\item
				$\MM(\vect w, \vect z) \geq \Phi(\vect w)$ for all $\vect w, \vect z \in \R^{Nn}$;

			\item
				$\MM(\vect z, \vect z) = \Phi(\vect z)$ for all $\vect z \in \R^{Nn}$.
			\end{enumerate}
			The \emph{Bregman proximal mapping} $\prox_\FBE: \R^{Nn} \rightrightarrows \R^{Nn}$ and its value function the \emph{Bregman Moreau envelope} $\FBE: \R^{nN} \to \R$ (recall the definitions in \eqref{eq:prox_def}) are then defined as
			\begin{equation}\label{eq:Moreau_envelope}
				\prox_\FBE(\vect z) \coloneqq \argmin_{\vect w \in \R^{Nn}} \MM(\vect w, \vect z), \quad
				\FBE(\vect z) \coloneqq \min_{\vect w \in \R^{Nn}} \MM(\vect w, \vect z).
			\end{equation}
			The corresponding \emph{forward-backward residual} is defined as $\Res = \id - \prox_\FBE$.
			The envelope $\FBE$ is considered as the Lyapunov function in our convergence studies.

			We proceed with the following fact that is the key to our convergence analysis.
			First, to facilitate our subsequent analysis, we introduce the matrix $U_i\in \R^{Nn\times n}$ whose $i$-th block rows form an identity matrix. For a given vector $y \in \R^{n}$, the action of $U_i$ can be expressed as:
			\begin{equation}\label{eq:U}
				U_i y
			=
				{(0,\dots, 0, \!\overbracket[0.5pt]{\,y\,}^{\mathclap{i\text{-th block}}}\!, 0, \dots, 0)}.
			\end{equation}
			The following fact demonstrates that block-coordinate updates result in descent on the Bregman Moreau envelope $\FBE$ in \eqref{eq:Moreau_envelope}.
			As remarked above this is not necessarily the case for the cost function.

			\begin{fact}[\protect{Descent lemma \cite[Lem. 4.2]{latafat2021bregman}}]\label{lem:descent_lemma_pre}
				Suppose that \cref{ass:basic} holds, and let $\vect x = (x_1,\dots, x_N)$.
				Fix $\vect y \in \prox_\FBE(\vect x)$, and let $\mathcal{I} \subseteq [N]$ be a subset of indices.
				Consider the block-coordinate update
				\begin{equation}
					\textstyle
					\bm v = \bm x + \sum_{i\in \mathcal I}U_i U_i^\top (\bm y -  \bm x),
				\end{equation}
				where $U_i$ is as defined in \eqref{eq:U}.
				Then,
				\begin{equation}
					\FBE(\vect v) \leq \FBE(\vect x) - \DB_{\hat{H}}(\vect v, \vect x).
					\label{eq:descent_lemma}
				\end{equation}
			\end{fact}
			This observation is then utilized to establish that the limit points of the sequence \((z^k)_{k\in\N}\) correspond to stationary points of the function \(\varphi\), which, in the nonconvex setting, represents the necessary condition $0 \in \hat\partial\varphi(z^\star)$.

			In the following fact we present some of useful properties of the Bregman proximal mapping $\prox_\FBE$ and the Bregman distance $\DB_{\hat{H}}$, and expand on their relation to those defined in \cref{sec:Alg}.

			\begin{fact}[\emph{Bregman proximal mapping} {\cite[Lem. 3.1]{latafat2021bregman}}]\label{lem:equivalence}%
				Suppose that \cref{ass:basic} holds and let $\gamma_i \in (0,\tfrac{N}{L_i})$.
				Then, the following hold:
				\begin{factenumerate}[label=(\roman*),ref=\ref*{lem:equivalence}(\roman*),widest=2]
				\item\label{thm:Legendre:hatHh}%
					\(
						\DB_{\hat{H}}(\bm v,\bm z)
					\geq
						\sum_{i=1}^N(\tfrac{1}{\gamma_i}-\tfrac{L_{f_i}}{N})\DB_{h_i}(v_i,z_i)
					\),
					for  \( (\bm v, \bm z) \in \R^{Nn} \times \R^{Nn} \).
				\item\label{thm:Legendre:rangeT}%
					\(
						\prox_\Phi^{\hat{H}}(\bm z)
					=
						\{(v,\cdots,v):v\in \Ts(\sum_{i=1}^N\nabla\hat{h}_i(z_i))\}
					\),
					where \(\Ts\) is as in \eqref{eq:t}, is a nonempty and compact subset of $\Delta$.
				\end{factenumerate}
			\end{fact}

			When $\bm z\in \Delta$, one has a lower-dimensional representation of the Bregman Moreau envelope $\FBE$, provided in the following corollary of \cref{lem:equivalence}.

			\begin{corollary}[lower-dimensional representations]\label{cor:lowerDim}%
				Let \cref{ass:basic} hold and let $\gamma_i \in (0,\tfrac{N}{L_i})$.
				Then, with the Bregman Moreau operator and the envelope associated with $\varphi$ in \eqref{eq:problem_formulation} given by
				\begin{align*}
					\Tt(z)
				={} &
					\argmin_{w\in\R^n}\left\{\varphi(w) + \DB_{\hat{h}}(w,z)\right\},
				\\
					\fbe(z)
				={} &
					\min_{w\in\R^n}\left\{\varphi(w) + \DB_{\hat{h}}(w,z)\right\}
					=
					\mathcal{L}(v,z),
					\;\text{with}\;
					v\in \Tt(z),
				\end{align*}
				it holds that $\Tt = \Ts \circ\nabla \hat{h} = \prox_{\varphi}^{\hat{h}}$, and  $\FBE(\bm z) = \fbe(z)$  for \(\bm z = (z, \ldots, z)\in \Delta\).
				Moreover,
				\begin{corenumerate}[label=(\roman*),ref=\ref*{cor:lowerDim}(\roman*),widest*=2]
				\item\label{thm:env:ld:hatHh}%
					\(
						\DB_{\hat{h}}(v,z)
					\geq
						\sum_{i=1}^N(\tfrac{1}{\gamma_i}-\tfrac{L_{f_i}}{N})\DB_{h_i}(v,z)
					\),
					for \( (v, z) \in \R^n \times \R^n \).

				\item\label{thm:env:ld:OC}%
					If \(v\in\prox_\varphi^{\hat{h}}(z)\), then
					\(
						\nabla\hat{h}(z)
						-
						\nabla\hat{h}(v)
					\in
						\hat\partial\varphi(v)
					\);
					the converse also holds when \(\varphi\) is convex.
				\end{corenumerate}
			\end{corollary}

			An important consequence of \cref{lem:equivalence} and its \cref{cor:lowerDim} is that the range of \(\prox_\FBE\) is a subset of the con\-sen\-sus set $\Delta$ (cf. \cref{thm:Legendre:rangeT}). Moreover, by \cref{thm:env:ld:OC}, to any fixed point of \(\prox_\FBE\) (or \(\prox_{\varphi}^{\hat{h}}\)) there corresponds a stationary point for the original problem.
			That is to say
			\begin{align*}
				z^\star \in \fix\prox_{\varphi}^{\hat{h}}
			\Leftrightarrow{}&
				0\in\res(z^\star)
				\Leftrightarrow
				z^\star \in \prox_{\varphi}^{\hat{h}}(z^\star)
			\stepcounter{equation}\tag{\theequation}
			\label{eq:fixed_point}
			\\
			\Leftrightarrow{}&
				(z^\star,\dots,z^\star) = \vect z^\star  \in \prox_\FBE(\vect z^\star)
			\\
			\Leftrightarrow{}&
				0 \in \Res(\bm z^\star)
				\Leftrightarrow
				0 \in \hat{\partial} \Phi(\vect z^\star)
			\\
			\Leftrightarrow{}&
				0\in \hat{\partial} \varphi(z^\star),
			\end{align*}
			where $\res$ is defined in \eqref{eq:FB_residual}.

		\subsection{Lifted Representation of the Algorithm}\label{sec:lifted}

			For the sake of clarity in presentation and without loss of generality, we consider the cyclic sweeping rule in the incremental loop where $i^\ell = \ell$ in step \ref{alg:1:rnd} (cf. \cref{rem:suffledCyclic} for shuffled cyclic sweeping rule).
			In this case, we adopt the following notation:
			\begin{equation}
				\tz[\ell]{k} \coloneqq (\tilde{z}^k_1,\tilde{z}^k_2,\dots,\tilde{z}^k_{\ell-1},\underbrace{u^k,\dots,u^k}_{N-\ell+1}), ~~\ell\in [N].
				\label{eq:tz_def}
			\end{equation}
			Using the defined operators in the previous subsection, the proposed \cref{alg:LBFinito} is outlined in the lifted space in \cref{tab:lifted}.

			\begin{algorithm}[t]
				\centering
				\caption{Representation of \cref{alg:LBFinito} in the lifted space}

				\begin{algorithmic}[1]
				\item[{\textbf{Initialize}}]
					\begin{tabular}[t]{@{}l@{}}
						\(\tz[N]{-1}\in\R^{Nn}\),~
						\(\beta \in (0,1)\),~
						\(\gamma_i\in(0,\nicefrac{N}{L_{i}})\), \(i\in[N]\)%
					\\
						maximum number of backtracks \(q_{\rm max} \in \N\cup\{\infty\}\){\tiny{(e.g. \(q_{\max}=2\))}},~
						\(K\in\N\)
					\end{tabular}

				\item[{\textbf{Repeat for}} $k=0,1,\dots,K$]
				\State \label{algLift:1:t1}%
					$\bm z^k \in \prox_\FBE(\tz[N]{k-1})$

				\State \label{algLift:1:t2}%
					$\bm v^k \in \prox_\FBE(\bm z^k)$ \hfill \Comment{(full update)}

				\State \label{algLift:1:dir}%
					choose $\bm d^k \in \Delta$

				\State \label{algLift:1:ls}%
					set $\tau_k = 1, q_k=0$
				\Comment{(linesearch)}
					\begin{enumerate}[label=\alph*:,ref=\ref*{algLift:1:ls}.\alph*,widest*=4]
					\item \label{algLift:2:ls:1}%
						\(\bm u^k=\tau_k \bm z^k+(1-\tau_k) \bm v^k+\tau_k{ \bm d^k}\)

					\item
						\(\bm y^k \in \prox_\FBE(\bm u^k)\)
					\Comment{(full update)}

					\item \label{algLift:1:ls:condition}%
						\textbf{if} \(\FBE(\bm u^k)\leq\FBE(\bm z^k)\)

						\quad go to step \ref{algLift:2:sts}

					\item
						\textbf{else if} \(q_k = q_{\max}\) \textbf{then}

						\quad
						\(\bm u^k = \bm v^k\), and go to step \ref{algLift:2:sts}

					\item
						\textbf{else}

						\quad
						\(\tau_k \gets \beta \tau_k\),~
						\(q_k \gets q_k + 1\),
						and go to step \ref{algLift:2:ls:1}
					\end{enumerate}

				\State \label{algLift:2:sts}%
					\( \tz[1]{k} = \bm u^k \)

				\For{ \(\ell = 1,\ldots, N\) }\label{algLift:forloop}%
				\Comment{(incremental loop)}
					\State \label{algLift:1:t3}%
						\(\tbz[\ell]{k} \coloneqq (\ttz[\ell]{k},\cdots,\ttz[\ell]{k}) \in \prox_\FBE(\tz[\ell]{k})\)

					\State \label{algLift:forloop:end}%
						\(\tz[\ell+1]{k} = \tz[\ell]{k} + U_{\ell}U_{\ell}^\top (\tbz[\ell]{k} - \tz[\ell]{k})\)
				\EndFor
				\end{algorithmic}
				\label{tab:lifted}
			\end{algorithm}

			In \cref{tab:lifted}, \(\bm z^k, \bm v^k\) belong to the consensus set \(\Delta\) owing to \cref{thm:Legendre:rangeT}.
			This along with the choice of \(\bm d^k\) implies the same for \(\bm u^k\), ensuring that the linesearch can be performed in the lower dimensional space (see \cref{rem:wellDefinedness}).
			In the following proposition, we highlight the equivalence of \cref{alg:LBFinito} and its lifted variant \cref{tab:lifted}.
			The proof is omitted as it follows directly from the above observations along with \cref{lem:equivalence} and \cref{cor:lowerDim}.

			\begin{proposition}\label{prop:equivalence}
				As long as the two algorithms are initialized with the same parameters, to any sequence $(z^k,v^k,d^k,u^k,y^k,\tilde z_{\ell}^k,\ell\in[N])_{k\in\N}$ generated by \cref{alg:LBFinito}, there correspond sequences $(\bm z^k = (z^k,\cdots,z^k))_{k\in\N}$, $(\bm v^k = (v^k,\cdots,v^k))_{k\in\N}$, $(\bm d^k=(d^k,\cdots,d^k))_{k\in\N}$, $(\bm u^k=(u^k,\cdots,u^k))_{k\in\N}$, $(\bm y^k=(y^k,\cdots,y^k))_{k\in\N}$, $(\tbz[\ell]{k}=(\ttz[\ell]{k},\cdots,\ttz[\ell]{k}))_{k\in\N}$, $\ell\in[N]$, generated by \cref{tab:lifted} (and vice versa).
			\end{proposition}

			Considering the updates at steps \ref{algLift:1:t1} and \ref{algLift:forloop:end}, the descent property established in \cref{lem:descent_lemma_pre} already hints as to why $\FBE$ is employed in the backtracking linesearch procedure.
			In the next remark, we expand on the well-definedness of this linesearch and discuss its relation to the one prescribed in \cref{alg:LBFinito}.
			These observations in the lifted space will help us in establishing convergence of the algorithm in \crefrange{sec:subseq}{sec:glb}.

			\begin{remark}[Well definedness of linesearch]\label{rem:wellDefinedness}
				The linesearch in step \ref{algLift:1:ls} is well defined, as it always terminates in a finite number of backtracks.
				In this process, the iterate $\bm u^k$ is interpolated between the candidate fast update $\bm z^k+ \bm d^k$, corresponding to \(\tau_k=1\), and the safeguard step $\bm v^k$ which is approached as \(\tau_k\searrow0\).
				Observe that \(\mathcal{L}(v^k,z^k)=\fbe(z^k)=\FBE(\bm z^k)\), as it follows from \cref{cor:lowerDim}, and that similarly \(\mathcal{L}(y^k,u^k)=\fbe(u^k)=\FBE(\bm u^k)\).
				Due to \cref{lem:descent_lemma_pre} and the continuity of the function $\FBE$, as long as $\bm z^k \neq \bm v^k$, the inequality $\FBE(\bm v^k) < \FBE(\bm z^k)$ holds, hence the inequality is satisfied for $\tau_k$ small enough.
				In practice, it is beneficial to limit the number of backtracks using a maximum value $q_{\max}$, especially at initial iterations.
				Regardless, as it will be shown in \cref{thm:superlinear}, under mild assumptions at the limit point eventually the iterate enters a region where backtracks will never be invoked.
			\end{remark}

			The updates in the incremental loop are referred to as block-coordinate updates, since at step \ref{algLift:forloop:end} only one block of $\tz[\ell]{k}$ is updated at each incremental loop iteration (equivalently as in step \ref{alg:forloop:end} of \cref{alg:LBFinito} due to \cref{thm:Legendre:rangeT}).
			Note that at each block update, due to choosing $\ell$ without replacement, the previous block value is known and equal to $u^k$, as depicted in \eqref{eq:tz_def}.
			Consequently, the update in step \ref{algLift:forloop:end} (equivalently in step \ref{alg:forloop:end} of \cref{alg:LBFinito}) can be accomplished by replacing $u^k$ with $\ttz[\ell]{k}$, thereby requiring a memory of $\mathcal{O}(n)$ instead of $\mathcal{O}(nN)$.
			However, opting for a memory of $\mathcal{O}(nN)$ and saving $\tilde{z}^k_{j}$ for all $j\in[N]$ allows the algorithm to adopt the randomized sweeping rule with replacement as well (cf. \cref{rem:incremental}).
			Additionally, the block updates in step \ref{algLift:1:t3} are computationally inexpensive, since if implemented by steps \ref{alg:1:t3} and \ref{alg:forloop:end} of \cref{alg:LBFinito} (using \cref{thm:Legendre:rangeT}), the gradient of only one function $f_{\ell}$ is computed in each incremental iteration.

			\begin{remark}[shuffled cyclic (randomized without replacement) sweeping rule]\label{rem:suffledCyclic}%
				It is evident from step \ref{algLift:forloop:end} that the incremental loop can be easily modified to accommodate the shuffled cyclic (randomized without replacement) sweeping rule, as was commented for \cref{alg:LBFinito}.
				To achieve this, we can consider the following update in place of step \ref{algLift:forloop:end}:
				\[
					\tz[\ell+1]{k} = \tz[\ell]{k} + U_{i^\ell}U_{i^\ell}^\top (\tilde{\bm z}_\ell^k - \tz[\ell]{k}),
				\]
				where $i^\ell \in [N]$ is randomly chosen without replacement.
			\end{remark}

		\subsection{Global and Subsequential Convergence}\label{sec:subseq}

			SPIRAL is globally (as opposed to locally) convergent whenever \cref{ass:basic} holds, without any additional assumption on the convexity of nonsmooth term $g$.
			Thanks to the proposed linesearch in step \ref{alg:1:ls}, global convergence is also guaranteed with any direction $d^k$ derived in step \ref{alg:1:dir}, although ultimately a fast convergence rate is only achieved by employing an educated direction and under assumptions at the limit point.
			Motivated by \cref{lem:descent_lemma_pre} and \cref{cor:lowerDim} and consistent with previous studies such as \cite{teboulle2018simplified,latafat2021bregman}, the Bregman Moreau envelope $\FBE$ in \eqref{eq:Moreau_envelope} is employed as the Lyapunov function which reduces to $\mathcal{L}(v^k,z^k)$ in the original space. This function has nice properties which enables us to study the global and subsequential convergence of SPIRAL, in the next theorem:

			\begin{theorem}{(Global and subsequential convergence)}\label{lem:descent_lemma}\label{thm:subsequential}
				Suppose that \cref{ass:basic} holds.
				The following holds for the sequence $(z^k)_{k\in\N}$ generated by \cref{alg:LBFinito}:
				\begin{thmenumerate}[label=(\roman*),ref=\ref*{thm:subsequential}(\roman*),widest*=3]
				\item \label{lem:descent_lemma:1}%
					$\mathcal{L}(v^{k+1},z^{k+1}) \leq \mathcal{L}(v^{k},z^{k}) - \sum_{i=1}^N \DB_{\hat{h_i}}(z^{k+1}, \tilde{z}^k_i)$ for $k \in \mathbb{N}$, with $\hat{h_i} = \frac{1}{\gamma_i}h_i - \frac{1}{N} f_i$;

				\item \label{lem:descent_lemma:res}%
					$(\DB_{\hat{h_i}}(z^{k+1}, \tilde{z}^k_i))_{k\in\N} \to 0, i \in [N]$;

				\item \label{lem:descent_lemma:3}%
					$(\mathcal{L}(v^{k},z^{k}))_{k\in\N}$ and $(\varphi(z^k))_{k\in\N}$ converge to a value $\varphi_\star$ where $\varphi(z^0) \geq \varphi_\star \geq  \inf \varphi$;

				\item \label{lem:descent_lemma:4}%
					$\varphi$ equals $\varphi_\star$ on all the cluster points;

				\item \label{lem:descent_lemma_stationary}%
					all the cluster points are fixed points for $\prox_\fbe$, and are in particular stationary for $\varphi$;

				\item \label{lem:descent_lemma:bounded}%
					if \(\varphi\) is level bounded, then $(z^k)_{k\in\N}$, $(\ttz[i]{k})_{k\in\N}$, for $i\in [N]$ are bounded.
				\end{thmenumerate}
			\end{theorem}
			\begin{proof}
				\ref{lem:descent_lemma:1}:
					The block-coordinate interpretation of vectors \(\tz[\ell]{k}\) as shown in step \ref{algLift:forloop:end} of \cref{tab:lifted} along with \cref{lem:descent_lemma_pre} yields
					\begin{equation}
						\FBE(\bar{\bm z}^k_{\ell+1}) \leq \FBE(\bar{\bm z}^k_{\ell}) - \DB_{\hat{H}}(\bar{\bm z}^k_{\ell+1}, \bar{\bm z}^k_{\ell})  \leq \FBE(\bar{\bm z}^k_{\ell}) \quad \text{for } \ell=1,\dots,N-1.
					\end{equation}
					By unrolling the inequality above we have
					\begin{equation}
						\FBE(\bar{\bm z}^k_N)
					\leq
						\FBE(\bar{\bm z}^k_{1})
					\leq
						\FBE(\vect u^k)
					\leq
						\FBE(\vect z^k),
						\label{eq:des_lem_1}
					\end{equation}
					where the second inequality uses \cref{lem:descent_lemma_pre} and the last one is ensured by the linesearch condition in step \ref{alg:1:ls}.
					Moreover, in step \ref{alg:1:t1}, $z^{k+1} = \Ts(\sum_{i=1}^N \nabla \hat{h}_{i}(\ttz[i]{k}))$, or equivalently stated  by \cref{thm:Legendre:rangeT}  $\vect z^{k+1} = (z^{k+1},\ldots, z^{k+1})\in \prox_\FBE(\tz[N]{k})$.
					Therefore, using \cref{lem:descent_lemma_pre} yields
					\begin{equation}
						\FBE(\vect z^{k+1}) \leq \FBE(\bar{\bm z}^k_{N}) - \DB_{\hat{H}}(\vect z^{k+1}, \bar{\bm z}^k_{N}).
						\label{eq:des_lem_2}
					\end{equation}
					Summing up the two inequalities in \eqref{eq:des_lem_1} and \eqref{eq:des_lem_2} yields
					\begin{equation}
						\FBE(\vect z^{k+1}) \leq \FBE(\vect z^k) - \DB_{\hat{H}}(\vect z^{k+1}, \bar{\bm z}^k_N).
						\label{eq:layapanov_decrease}
					\end{equation}
					Noting that $\vect z^k = (z^k,\dots,z^k)$, the above inequality may be written as (cf. \cref{cor:lowerDim})
					\begin{equation*}
						\textstyle \mathcal {L}(v^{k+1},z^{k+1}) \leq \mathcal {L}(v^{k},z^{k}) - \sum_{i=1}^N \DB_{\hat{h_i}}(z^{k+1}, \tilde{z}^k_i).
					\end{equation*}

				\ref{lem:descent_lemma:res}:
					By reordering the inequality in \eqref{eq:layapanov_decrease} and telescoping we have
					\begin{equation*}
						\textstyle \sum_{k=0}^T \DB_{\hat{H}}(\vect z^{k+1}, \bar{\bm z}^k_N) \leq \FBE(\vect z^0) - \FBE(\vect z^T) \leq \FBE(\vect z^0) - \inf \FBE < \infty.
					\end{equation*}
					The last two inequalities follow from the boundedness of $\FBE$ from below, in light of \cref{ass:basic} and \cref{prop:Moreau_properties_3}.
					The inequality above shows that the sum is finite and hence $(\DB_{\hat{H}}(\vect z^{k+1}, \bar{\bm z}^k_N))_{k\in\N} \to 0$.

				\ref{lem:descent_lemma:3}:
					The sequence $(\FBE(\vect z^k))_{k\in\N}=(\mathcal {L}(v^{k},z^{k}))_{k\in\N}$ is decreasing by \eqref{eq:layapanov_decrease} and since it is lower bounded, it should converge to a finite value $\varphi_\star$ with $\varphi(z^0) \geq \varphi_\star \geq \inf \FBE = \inf {\Phi} =  \inf \varphi$.
					Moreover, from \eqref{eq:des_lem_2} and \eqref{eq:des_lem_1}:
					\[
						\FBE(\vect z^{k+1}) + \DB_{\hat{H}}(\vect z^{k+1}, \bar{\bm z}^k_{N})
					\leq
						\FBE(\bar{\bm z}^k_{N})
					\leq
						\FBE(\vect z^{k}).
					\]
					Since \(\DB_{\hat{H}}(\vect z^{k+1}, \bar{\bm z}^k_{N})\) vanishes, see \cref{lem:descent_lemma:res}, \(\FBE(\bar{\bm z}^k_{N}) \to \varphi_\star\), which in turn implies through the identity
					\(
						\FBE(\bar{\bm z}^k_{N})
					=
						\Phi(\bm z^{k+1}) + \DB_{\hat{H}}(\vect z^{k+1}, \bar{\bm z}^k_{N})
					\),
					that $(\Phi(\vect z^k))=(\varphi(z^k))_{k\in\N} \to \varphi_\star$.

				\ref{lem:descent_lemma:4}:
					Take a subsequence $(\vect z^k_{k\in K} \to \vect z^\star$ with $K \subseteq \N$.
					We have:
					\begin{equation}
						\varphi_\star
					\xleftarrow[k\in K]{\text{\ref{lem:descent_lemma:3}}}
						\FBE(\vect z^k) - \DB_{\hat{H}}(\vect z^\star, \vect z^k)
					\stackrel{\text{\ref{prop:Moreau_properties_1}}}{\leq}
						\Phi(\vect z^\star)
					\stackrel{\text{lsc}}{\leq}
						\lim_{k \in K} \inf \Phi(\vect z^k)
					\stackrel{\text{\ref{lem:descent_lemma:3}}}{=}
						\varphi_\star.
					\end{equation}

				\ref{lem:descent_lemma_stationary}:
					Let $K \subseteq \N$ denote an infinite subsequence such that $(\vect z^k_{k \in K} \to \vect z^\star$.
					It follows from \cref{lem:descent_lemma:res} along with \cite[Thm. 2.4]{solodov2000inexact} that $(\tz[N]{k})_{k \in K} \to \vect z^\star$.
					With $\vect z^{k+1} \in \prox_\FBE(\bar{\bm z}^k_N)$ and the osc property of $\prox_\FBE$ (see \cref{prop:prox_properties}), it follows that $\vect z^\star \in \prox_\FBE(\vect z^\star)$ implying stationarity of the limit points as shown in \eqref{eq:fixed_point}.

					\ref{lem:descent_lemma:bounded}: Level boundedness of $\varphi$ implies that of $\FBE$.
					It then  follows from \eqref{eq:des_lem_2} and \eqref{eq:des_lem_1} that \((\tz[N]{k})_{k\in\N}\) is
					contained in \(\{\bm w:\FBE(\bm w) \leq \FBE({\tz[N]{0}})\}\), which is a bounded set. Boundedness of $(\bm z^k=(z^k,\ldots, z^k))_{k\in\N}$ follows from that of \((\tz[N]{k})_{k\in\N}\), local boundedness of the proximal mapping (see \cref{prop:prox_properties}), and $\vect z^{k+1} \in \prox_\FBE(\bar{\bm z}^k_N)$.
			\end{proof}

		\subsection{Superlinear Convergence}\label{sec:spl}

			In this section, we aim to demonstrate that the proposed linesearch is \emph{smart}, particularly in identifying \emph{mature} directions. When a direction $d^k$ is deemed \emph{mature}, the candidate update $z^k+d^k$ will eventually be accepted without any backtracks, thereby enabling SPIRAL to exhibit superlinear convergence.

			We first introduce necessary additional assumptions and lemmas.
			Subsequently, we delve into the identification of \emph{mature} directions by SPIRAL and explore their relationship with quasi-Newton methods.
			Take the following assumptions:

			\begin{assumption}[superlinear convergence requirements]\label{ass:superlinear}%
				The following hold in problem \eqref{eq:problem_formulation}:
			    \begin{assenumerate}[label=A\arabic*,ref=\ref*{ass:superlinear}.A\arabic*]
				\item \label{ass:g:cvx}%
					\(g\) is convex;
				\item \label{ass:fh}%
					for $i\in[N]$, \(f_i,h_i\in\C^2\) with \(\nabla^2h_i\succ 0\).
				\end{assenumerate}
			\end{assumption}

			To achieve superlinear convergence, SPIRAL requires that the sequence $(z^k)_{k\in\N}$ converges to a strong local minimum $z^\star$ of the cost function, and that the envelope is twice (strictly) differentiable at $z^\star$.
			It is noteworthy that the strong local minimality (isolated local minima) is a standard requirement for asymptotic properties of quasi-Newton methods.
			However, works such as \cite{aragon2014local,themelis2019supermann} relax this requirement to address nonisolated local minima as well.
			As future work, their techniques can be investigated for our setting.

			To establish the aforementioned properties of the envelope, the following fact and lemma are presented, taking into account the additional \cref{ass:superlinear} in conjunction with \cref{ass:basic}.
			First, let us formally define the concept of strong local minimality:

			\begin{definition}[Strong local minimum]
				A point $z^\star$ is said to be the strong local minimum of $\phi$ if there exist a neighborhood $\mathcal N_{z^\star}$ of $z^\star$ and $c>0$ such that for all $z\in \mathcal N_{z^\star}$, $\phi(z) \geq \phi(z^\star) + \tfrac{c}2\|z-z^\star\|^2$.
			\end{definition}

			The following fact establishes an equivalence between strong local minima of function $\varphi$ and of its envelope $\fbe$.

			\begin{fact}[equivalence of strong local minima \protect{\cite[Thm. 3.7]{ahookhosh2021bregman}}]\label{lem:equiv:strloc}%
				Suppose that \cref{ass:basic,ass:superlinear} hold.
				Then, $z^\star\in \fix \prox_{\varphi}^{\hat{h}}$ is a strong local minimum of $\varphi$ if and only if it is a strong local minimum of $\fbe$.
			\end{fact}

			We remark that $\prox_{\varphi}^{\hat{h}}$ is single-valued whenever $g$ is convex, hence all the fixed points are guaranteed to be \emph{nondegenerate} in the sense of \cite[Def. 3.5]{ahookhosh2021bregman}.
			In order to achieve superlinear convergence, we assume $z^\star$ to be a strong local minimum of the cost $\varphi$, and that the envelope is twice (strictly) differentiable at this point. The subsequent lemma examines the second-order properties of the envelope to ensure the fulfillment of this requirement.

			\begin{lemma}[second order characterization]\label{lem:second:fbe}%
				Suppose that \cref{ass:basic,ass:superlinear} hold.
				Then, given $z^\star \in \fix \prox_\fbe$,  there exists a neighborhood of $z^\star$ where $\res$ is Lipschitz continuous, and $\fbe$ is continuously differentiable.
				If, in addition, $\Ts$ is (strictly) differentiable at \(\nabla\hat{h}(z^\star)\), then
				\begin{lemenumerate}[label=(\roman*),ref=\ref*{lem:second:fbe}(\roman*),widest*=2]
				\item \label{eq:res:Jac}%
					\(\res\) is (strictly) differentiable at $z^\star$ with
					\(
						\jac{\res}{(z^\star)} = \I - \jac{\Ts}{(\nabla \hat{h}(z^\star))}\nabla^2\hat{h}(z^\star);
					\)

				\item \label{eq:HessianEnv}%
					\(\fbe\) is twice (strictly) differentiable at $z^\star$ with symmetric Hessian
					\[
						\nabla^2\fbe(z^\star) = \nabla ^2\hat{h}(z^\star) \jac{\res}{(z^\star)}.
					\]
					In particular, if $z^\star$ is a strong local minimum of $\varphi$, then $\nabla^2 \fbe(z^\star)$ is symmetric positive definite and $\jac{\res}{(z^\star)}$ is invertible.
				\end{lemenumerate}
			\end{lemma}
			\begin{proof}
				Observe that $\Tt = \Ts \circ\nabla \hat{h} = \prox_{\varphi}^{\hat{h}}$ as shown in \cref{cor:lowerDim}.
				Given this characterization, the first claim follows directly from \cite[Thm. 3.10]{ahookhosh2021bregman}.

				\ref{eq:res:Jac}
					Since $\nabla \hat{h}\in \C^1$ (cf. \cref{ass:fh}) and $\Ts$ is (strictly) differentiable at $\nabla\hat{h}(z^\star)$, so is the composition \(\Ts\circ \nabla \hat{h}\), thus implying (strict) differentiability of $\res = \id - \Tt$.
					The Jacobian of the residual is obtained by the chain rule.

				\ref{eq:HessianEnv}:
					The claim follows from (strict) differentiability of $\Ts$ at $\nabla\hat{h}(z^\star)$.
					Moreover, $\nabla^2 \hat{h}\succ0$ owing to \cref{ass:basic:1,ass:fh} and $\gamma_i\in(0,\nicefrac{N}{L_i})$. Thus, $\jac{\res}{z^\star}$ is nonsingular when $z^\star$ is a strong local minimum of $\fbe$, or, equivalently, of $\varphi$ (cf. \cref{lem:equiv:strloc}).
			\end{proof}

			After establishing the desirable properties of the envelope in the vicinity of fixed point $z^\star\in\fix\prox_\fbe$, we now proceed to characterize the quality of directions $d^k$ in step \ref{alg:1:dir} of \cref{alg:LBFinito} through introducing the notion of \emph{superlinear directions} which was introduced in the seminal work \cite{facchinei2003finite}.
			We also refer the reader to the works such as \cite{themelis2019supermann,ahookhosh2021bregman} for further discussion and extensions.

			\begin{definition}[superlinear directions]\label{def:superlinearDirection}%
				Relative to a sequence $(z^k)_{k\in\N}$ that converges to a point $z^\star\in\R^n$, we say that $(d^k)_{k\in\N}\subseteq \R^n$ is a sequence of superlinear directions, if
				\[
					\lim_{k\to\infty}{\frac{\|z^k+d^k - z^\star\|}{\|z^k- z^\star\|}} = 0.
				\]
			\end{definition}

			The following theorem provides asymptotic guarantees for the superlinear convergence of SPIRAL.
			As remarked before, when the directions satisfy \cref{def:superlinearDirection}, the backtracking linesearch will eventually never be triggered, thus substantially improving the performance of the algorithm.
			The theorem requires local assumptions such as strict differentiability of $\Ts$ in \eqref{eq:t} at $\nabla \hat{h}(z^\star)$, where $z^\star$ is the limit point of $(z^k)_{k\in\N}$.
			Auxiliary results for controlling the terms $\|\ttz[\ell]{k} - u^k\|$ for $\ell\in[N]$ that appear due to the (shuffled) cyclic updates are postponed to \cref{lem:Tri} in Appendix \ref{sec:omitted_lemmas}.

			\begin{theorem}[superlinear convergence]\label{thm:superlinear}
				Consider the sequence $(z^k)_{k\in\N}$ generated by \cref{alg:LBFinito}, and additionally to \cref{ass:basic,ass:superlinear}, suppose the following are satisfied:%
			    \begin{thmenumerate}[label=A\arabic*,ref=\ref*{thm:superlinear}.A\arabic*]
				\item \label{thm:ass:strloc}%
					$(z^k)_{k\in\N}$ converges to a strong local minimum $z^\star$ of $\varphi$;

				\item
					the directions $d^k$ are superlinear relative to $(z^k)_{k\in\N}$ (cf. \cref{def:superlinearDirection});

				\item \label{thm:ass:tdiff}%
					$\Ts$  (defined in \eqref{eq:t}) is strictly differentiable at $\nabla \hat{h}(z^\star)$.
				\end{thmenumerate}
				Then, asymptotically the linesearch in step \ref{alg:1:ls} will be accepted with $\tau = 1$, and $(z^k)_{k\in\N}$ converges to $z^\star$ at superlinear rate.%
			\end{theorem}
			\begin{proof}
				Let $u_0^k \coloneqq z^k + d^k$.
				Due to the superlinearity of the directions $d^k$,
				\begin{equation}
					\lim_{k\to \infty}\tfrac{\|u_0^k-z^\star\|}{\|z^k-z^\star\|}
					= 0.
					\label{eq:suplim:u0}
				\end{equation}
				We start by showing that close enough to the limit point the linesearch condition would always be satisfied with $\tau = 1$.
				It follows from \cref{lem:descent_lemma_stationary} that \(z^\star \in \fix \prox_\fbe\) and from \cref{lem:equiv:strloc} that $z^\star$ is also a strong local minimum of $\fbe$.
				Hence, $G^\star = \nabla^2 \fbe(z^\star)$ is symmetric positive definite by \cref{eq:HessianEnv}.
				Let
				\[
					\varepsilon_k \coloneqq \tfrac{\fbe(u_0^k)-\varphi_{\star}}{\fbe(z^k)-\varphi_{\star}}.
				\]
				Since $G^\star \coloneqq \nabla^2\fbe(z^\star)\succ0$, a second-order expansion of \(\fbe\) at $z^\star$ yields
				\begin{align*}
					\lim_{k\to \infty}\varepsilon_k
				={} &
					\lim_{k\to \infty}\frac{\tfrac12\langle G^\star(u_0^k-z^\star), u_0^k - z^\star \rangle + \mathcal{O}\big(\|u_0^k-z^\star\|^2\big)}{\tfrac12\langle G^\star(z^k-z^\star), z^k - z^\star \rangle + \mathcal{O}\big(\|z^k-z^\star\|^2\big)}
				\\
				\leq{}&
					\lim_{k\to \infty}\frac{\|G^\star\|\|u_0^k-z^\star\|^2 + \mathcal{O}\big(\|u_0^k-z^\star\|^2\big)}{\lambda_{\rm min}(G^\star)\|z^k-z^\star\|^2 + \mathcal{O}\big(\|z^k-z^\star\|^2\big)}
				\\
				={}&
					\lim_{k\to \infty}\frac{\|G^\star\|\tfrac{\|u_0^k-z^\star\|^2}{\|z^k-z^\star\|^2} + \tfrac{\mathcal{O}\big(\|u_0^k-z^\star\|^2\big)}{\|z^k-z^\star\|^2}}{\lambda_{\rm min}(G^\star) + \tfrac{\mathcal{O}\big(\|z^k-z^\star\|^2\big)}{\|z^k-z^\star\|^2}}
				\stackrel{\text{\eqref{eq:suplim:u0}}}{\leq}
					0.
				\end{align*}
				In particular, there exists $k_0\in\N$ such that $\varepsilon_k \leq 1 ~\forall k \geq k_0$.
				Moreover, since $z^k$ converges to $z^\star$, it follows from \cref{lem:descent_lemma:1} and \cref{cor:lowerDim} that
				\(%
					\fbe(z^k) \geq \fbe(z^\star).
				\) %
				Consequently, using the definition of $\varepsilon_k$ above with $\varphi_{\star} = \fbe(z^\star)$ due to \cref{prop:Moreau_properties_3},
				\[
					\fbe(u_0^{k}) - \fbe(z^k) = -(1-\varepsilon_k)
					\left(%
						\fbe(z^k) - \fbe(z^\star)
					\right)
					\leq 0 \quad \forall k \geq k_0.
				\]
				Therefore, the unit stepsize would always be accepted in step \ref{alg:1:ls}.
				Now, with unit stepsize $u_0^k=u^k$, as in step \ref{alg:1:ls:1}, and this implies through \eqref{eq:suplim:u0} that
				\begin{equation}
					\lim_{k\to \infty}\frac{\|u^k - z^\star\|}{\|z^k-z^\star\|} = 0.
					\label{eq:super_linear_lim}
				\end{equation}
				On the other hand
				\begin{align*}
					\|\tz[N]{k} - \vect u^k\|
				\stackrel{\text{\eqref{eq:tz_def}}}{\leq} &
					\textstyle \sum_{i=1}^N \|\ttz[i]{k} - u^k\|
				\stackrel{\eqref{eq:xu}}{\leq}
					\eta \|\ttz[1]{k} - u^k\|
				\\
					\text{\footnotesize triangular inequality }
				\leq{} &
					\eta \|\ttz[1]{k} - z^\star\| + \eta \|u^k - z^\star\|
				\\
					\text{\footnotesize step \ref{algLift:1:t3} of \cref{tab:lifted}}
				={} &
					\textstyle \frac{\eta}{\sqrt{N}} \|\tbz[1]{k} - \bm z^\star\| + \eta \|u^k - z^\star\|
				\\
				\leq{} &
					\eta (\bar L + 1) \|u^k - z^\star\|,
				\stepcounter{equation}\tag{\theequation}
				\label{eq:zuzstar}
				\end{align*}
				where $\eta = \sum_{i=1}^N c_i$ and the last inequality follows from local Lipschitz continuity of $\prox_\FBE$ and step \ref{algLift:2:sts} of \cref{tab:lifted}.
				Further exploiting local Lipschitz continuity of the proximal mapping
				\begin{align*}
					\|z^{k+1} - z^\star\|
				={}&
					\tfrac{1}{\sqrt{N}}\|\vect z^{k+1} - \vect z^\star\|
				\\
					\text{\footnotesize Lip. cont. of $\prox_\FBE$ and step \ref{algLift:1:t1} of \cref{tab:lifted}}
				\leq{} &
					\tfrac{\bar L}{\sqrt{N}} \|\tz[N]{k} - \vect z^\star\|
				\\
					\text{\footnotesize triangular inequality }
				\leq{}&
					\tfrac{\bar L}{\sqrt{N}} \|\tz[N]{k} - \vect u^k\| + \tfrac{\bar L}{\sqrt{N}}\|\vect u^k - \vect z^\star\|
				\\
					\text{\footnotesize \eqref{eq:zuzstar}}
				\leq{} &
					\alpha \|u^k - z^\star\|
					\label{eq:super_linear_rel}
				\end{align*}
				where $\alpha = \tfrac{\eta}{\sqrt{N}} \bar L(\bar L+1) + \bar L$.
				Hence, combined with \eqref{eq:super_linear_lim}
				\[
					\frac{\|z^{k+1} - z^\star\|}{\|z^k-z^\star\|} \leq \alpha \frac{\|u^k - z^\star\|}{\|z^k-z^\star\|} \to 0,
				\]
				establishing the claimed superlinear convergence.
			\end{proof}

			A well-known condition for analyzing quasi-Newton methods is the celebrated Dennis-Mor\'e condition \cite{dennis1974characterization,dennis1977quasi}, which characterizes the quality of the directions as follows:
			\begin{equation}
			\lim_{k\to \infty}{\frac{\|\res(z^k)+ \jac{\res}{(z^\star)}d^k\|}{\|d^k\|}}=0.
			\label{eq:dennis_more}
			\end{equation}
			This classical condition, in conjunction with \cref{ass:basic,ass:superlinear,thm:ass:strloc,thm:ass:tdiff}, leads to the emergence of superlinear directions \cite[Thm. 5.13]{ahookhosh2021bregman}.

			Note that the directions computed by Broyden updates, as one of the quasi-Newton methods, provably satisfy the Dennis-Mor\'e condition stated in \eqref{eq:dennis_more} (refer to \cite[Thm. 5.11]{themelis2018forward} and \cite[Thm. VI.8]{themelis2019supermann}).
			In order to achieve this, Broyden updates require the aforementioned regularity conditions on $\res$ at $z^\star$, as well as the boundedness of low-rank updates $H^k$ in \eqref{eq:update_rule_of_d}.
			It is important to mention that, although it is not formally established that L-BFGS satisfies the Dennis-Mor\'e condition, L-BFGS performs better than Broyden updates in practical scenarios.
			The theoretical examination of the Dennis-Mor\'e condition using L-BFGS updates is considered as a future research direction.

			It is noteworthy that the boundedness of low-rank updates $H^k$---similarly existence of the bound $\|d^k\|\leq D\|z^k-v^k\|$ with some finite $D \geq \sup_{k\in\N}{\|H^k\|}$ due to \eqref{eq:update_rule_of_d}, as it will be required in \cref{thm:global}---is a common assumption in the analysis of quasi-Newton methods (see, e.g., \cite[Ass. 2]{themelis2019supermann} and \cite[Thm. 5.7-A3 and Thm. 5.8-A3]{ahookhosh2021bregman}), and is guaranteed by employing safeguards in practice.

			\begin{remark}[practical considerations]\label{rem:dbound}
				The condition $\|d^k\|\leq D\|z^k-v^k\|$ is mild in practice since, as a safeguard here, in the case of failure in meeting the inequality, the directions may be scaled by $d^k \gets D\tfrac{\|z^k-v^k\|}{\|d^k\|}d^k$ using a sufficiently large predefined scalar $D$.
				Moreover, failure in meeting the inequality does not deteriorate the global and subsequential convergence of SPIRAL, as long as $\|d^k\|$ is scaled whenever necessary, as demonstrated in \cref{thm:subsequential}.
				As a final remark, although failure is possible (especially in the initial iterations when $d^k$ is not \emph{mature}), it has not occurred in any of our simulations in \cref{sec:sim}.
			\end{remark}

			\begin{remark}\label{rem:globalization}
				We remark that global convergence results in \cref{thm:subsequential} is established for any choice of direction.
				Even though in \cref{alg:LBFinito} quasi-Newton directions based on the residual mapping were suggested (cf. \eqref{eq:update_rule_of_d}), any superlinear direction can be employed in the algorithm. As a result, our theory provides a direct globalization strategy for works that employ quasi-Newton direction with only local convergence guarantees.
				For instance, it globalizes the recent work \cite{mokhtari2018iqn} which studies smooth and strongly convex finite sum problems, and proposes an incremental quasi-Newton method with local convergence guarantees.
			\end{remark}

		\subsection{Sequential and Linear Convergence}\label{sec:glb}

			In accordance with \cref{thm:superlinear} presented in the previous subsection, in order to achieve superlinear convergence, SPIRAL requires to have a sequence $(z^k)_{k\in\N}$ that converges to a strong local minimum of the cost function $\varphi$.
			The subsequent theorem establishes conditions under which the entire sequence $(z^k)_{k\in\N}$ converges to a stationary point with a linear convergence rate. To accomplish this, an additional assumption is required, namely, the Kurdyka-\L ojasiewicz (KL) property of the full cost function \cite{kurdyka1998gradients}.
			It is worth noting that $\varphi$ possesses the KL property for a wide range of problems, including situations where $f_i$ and $g$ are semialgebraic functions, which are commonly encountered in various applications (refer to \cite{bolte2007lojasiewicz,bolte2007clarke} for further elaboration).
			The formal statement of the KL property is as follows:

			\begin{definition}{(KL property with exponent $\theta$)}\label{def:KL_def}
				A proper lsc function $\phi: \R^n \to \Rinf$ has the Kurdyka-Lojasiewicz (KL) property with exponent $\theta\in (0,1)$ if for every $z^\star \in \dom \partial \phi$ there exist constants $\eta, \epsilon, \rho>0$ such that
				\begin{equation}
					\psi' (\phi(z)-\phi(z^\star))\dist(0,\partial \phi(z)) \geq 1, \quad \psi(s) \coloneqq \rho s^{1-\theta},
					\label{eq:KL_def}
				\end{equation}
				for all $z$ such that
				\(
					\|z-z^\star\| \leq \epsilon
				\)
				and
				\(
					\phi(z^\star) < \phi(z) < \phi(z^\star) + \eta
				\).
			\end{definition}

			In the following result, we present convergence guarantees for \emph{sequential} convergence under the KL assumption on the cost function. The proof follows standard techniques found in the literature, drawing inspiration from works such as \cite{Themelis2019Acceleration,ahookhosh2021bregman,latafat2021bregman}. For the sake of completeness, we provide the proof in Appendix \ref{sec:auxiliary:glb}.
			The primary challenge in establishing the \emph{sequential} convergence of \cref{alg:LBFinito} lies in demonstrating the bound $\|z^{k+1}- z^k\| \leq C\|\bm z^k - \tz[N]{k-1}\|$ for a positive constant $C$, as required by \cref{thm:global}.
			This bound is established by \cref{lem:Tri:2} presented in Appendix \ref{sec:omitted_lemmas}.

			\begin{theorem}[Sequential and linear convergence]\label{thm:global}%
				Additionally to \cref{ass:basic,ass:superlinear}, suppose the following is satisfied:
				\begin{thmenumerate}[label=A\arabic*,ref=\ref*{thm:global}.A\arabic*]
				\item
					$\varphi$ is level bounded;

				\item
					$\varphi$ has the KL property (cf. \cref{def:KL_def}) with exponent $\theta \in (0,1)$;

				\item
					the directions $d^k$ in step \ref{alg:1:dir} satisfy $\|d^k\|\leq D\|z^k-v^k\|$ for some \(D\geq 0\).
				\end{thmenumerate}
				Then, \((z^k)_{k\in\N}\) converges to a stationary point \(z^\star\) for \(\varphi\).
				Moreover, if the KL function has the exponent parameter in the range \(\theta\in(0,\nicefrac12]\), then \((z^k)_{k\in\N}\) and \((\varphi(z^k))_{k\in\N}\) converge at R-linear rate.
			\end{theorem}
			\begin{proof}
				Refer to Appendix \ref{aux:thm:global}.
			\end{proof}

	\section{Numerical Experiments}\label{sec:sim}

		In this section, we evaluate the proposed algorithm, SPIRAL, for both convex and nonconvex problems, considering cost functions with and without Lipschitz continuous gradients. We examine two versions of SPIRAL: 1) SPIRAL, which follows \cref{alg:LBFinito}, and 2) adaSPIRAL, an adaptive version with additional steps as outlined in \cref{tab:adaptive}.
		We compare SPIRAL against proxSARAH \cite{pham2020proxsarah}, proxSVRG \cite{reddi2016proximal}, proxSGD \cite{ghadimi2016mini}, proxSAGA \cite{reddi2016proximal}, Finito/MISO \cite{defazio2014finito,mairal2015incremental}, and low-memory Finito/MISO \cite[Alg. 2]{latafat2021bregman}.
		For the convex $\ell_1$ regularized least squares problem, we compare against Finito/MISO \cite{defazio2014finito}.
		For the nonconvex nonnegative principal component analysis problem, we compare against \cite{latafat2021block}, which addresses Finito/MISO in the general nonsmooth nonconvex case.
		Additionally, we compare SPIRAL against SMD and the Bregman Finito/MISO method \cite{latafat2021bregman} for the phase retrieval problem, where the cost function lacks a Lipschitz continuous gradient.
		The databases used in the evaluations are from LIBSVM \cite{chang2011libsvm}.
		To assess the performance of the algorithms, we employ the suboptimality criterion
		\begin{equation}
		    \mathcal{D}(z^k) \coloneqq \|z^k-v^k\|
		    \label{eq:optimality_criteria}
		\end{equation}
		with $v^k \in \Tt(z^k)$, since
		\begin{align*}
			\dist(0,\hat\partial \varphi(v^k))
		\leq{} &
			\inf_{v^k \in \Tt(z^k)} \|\sum_{i=1}^N \nabla \hat{h}_i(z^k) - \nabla \hat{h}_i(v^k)\|
		\\
		\leq{} &
			\sum_{i=1}^N \|\nabla \hat{h}_i(z^k) - \nabla \hat{h}_i(v^k)\|
			\leq c \|z^k-v^k\|,
		\end{align*}
		where the first inequality holds by \cref{thm:env:ld:OC} in \cref{sec:reform}, and $c > 0$ is some constant due to local Lipschitz continuity of $\hat{h}_i$, and the fact that $z^k$ and $v^k$ remain bounded.

		For all the algorithms in the comparisons, the stepsizes are set according to their theoretical convergence studies. Refer to \cite[Thm. 8]{pham2020proxsarah} for proxSARAH, \cite[Thm. 1]{reddi2016proximal} for proxSVRG, \cite[Thm. 3]{reddi2016proximal} for nonconvex proxSAGA, and \cite{defazio2014saga} for convex proxSAGA.
		For proxSGD the diminishing stepsize $\gamma^t = \frac{\gamma_0}{1+ \tilde{\gamma} t}$ is considered according to \cite{ghadimi2013stochastic} with $t$ as the epoch counter, $\gamma_0=0.1$, and $\tilde{\gamma}=0.5$. For (Bregman) Finito/MISO and also all the SPIRAL versions the stepsizes are set $\gamma_i=\frac{\alpha N}{L_i}, \;i \in [N]$ with $\alpha=0.999$.
		For adaSPIRAL, the stepsizes are all initialized by $\kappa \times \max_{i\in[N]}\{\frac{N}{L_i}\}$, with a grid search for $\kappa\in\{5, 10, 50, 100\}$ for each plot.
		Furthermore, the quasi-Newton directions in step \ref{alg:1:dir} are computed using \eqref{eq:update_rule_of_d}, where $H^k$ is updated using the L-BFGS method with a memory size of $5$.
		The maximum number of backtracks $q_{\max}$ is also set equal to $5$.
		It should be mentioned that while directions in step \ref{alg:1:dir} can be computed using any quasi-Newton method (e.g. Broyden updates as discussed in \cref{sec:spl}), L-BFGS yields superior numerical results.
		Finally, we refer to epochs by counting the total number of individual gradient evaluations divided by $N$, including those involved in the linesearches.

		\subsection{Adaptive variant}

			Using the global smoothness constants may lead to conservative stepsizes.
			\cref{tab:adaptive} expands steps \ref{alg:1:t1}, \ref{alg:1:t2}, \ref{alg:1:ls}, and \ref{alg:1:t3} of \cref{alg:LBFinito} in order to clarify the order of operations with the added backtracking linesearches that ensure the fundamental descent property in \cref{fact:smoothness}.
			This adaptation allows SPIRAL to estimate (relative) smoothness moduli locally---as opposed to using global estimates---resulting in larger stepsizes.
			The reader is referred to Appendix \ref{sec:auxiliary:adaptive} for further explanations about the memory-efficient implementation of the adaptive variant.

			\begin{table}[t]
				\centering
				\caption[]{%
					Adaptive stepsize selection for SPIRAL. Let $\alpha \in (0,1)$, $\sigma \in (0,1)$, and $\beta$ as in \cref{alg:LBFinito}.
					Initialize the stepsizes $\gamma_i$ for $i\in[N]$. For the sake of readability, define
					\(
						p_i(y,x)
					\coloneqq
						f_i(y) - f_i(x) - \langle \nabla f_i(x),y-x\rangle
					\).
				}

			\hrule
			\hrule
			\begin{algorithmic}[1]
				\State \label{state:tab1}%
				\begin{enumerate}[label=\alph*:,ref=\ref*{state:tab1}.\alph*,widest*=4]
				\item \label{tab:1:1}%
					\(z^k \in \Ts(s^k)\)

				\item \label{tab:1:2}%
					\textbf{if}
					\(
						\sum_{i=1}^N p_i(z^k,\ttz[i]{k-1})
						>
						\sum_{i=1}^{N}\tfrac{\alpha N}{\gamma_i}\DB_{h_i}(z^k,\ttz[i]{k-1})
					\)
					\textbf{then}

					\quad
					\(\gamma_i \gets \sigma \gamma_i\) for all \(i \in [N]\),
					update $s^k$, and go to step \ref{tab:1:1}%
				\end{enumerate}
			\end{algorithmic}

			\hrule

			\begin{algorithmic}[1]
			\addtocounter{ALG@line}{2}
				\State \label{state:tab2}%
				\begin{enumerate}[label=\alph*:,ref=\ref*{state:tab2}.\alph*,widest*=4]
				\item \label{tab:2:1}%
					\(v^k \in \Ts(\bar s^k)\)

				\item
					\textbf{if}
					\(
						\sum_{i=1}^N  p_i(v^k,z^k)
						>
						\sum_{i=1}^N \tfrac{N}{\gamma_i}\DB_{h_i}(v^k,z^k)
					\)
					\textbf{then}

					\quad
					\(\gamma_i \gets \sigma \gamma_i\) for all \(i \in [N]\),
					update $\bar s^k$, and go to step \ref{tab:2:1}%
				\end{enumerate}
			\end{algorithmic}

			\hrule

			\begin{algorithmic}[1]
			\addtocounter{ALG@line}{3}
			\State \label{tab:3:1}%
				choose $d^k\in\R^n$ at $z^k$ {\footnotesize(e.g. based on a quasi-Newton method for solving $\res(z)=0$)}

			\State \label{tab:3}%
				set $\tau_k = 1, q_k=0$
				\begin{enumerate}[label=\alph*:,ref=\ref*{tab:3}.\alph*,widest*=4]
				\item \label{tab:3:2}%
					\(u^k=z^k+(1-\tau_k)(v^k-z^k)+\tau_k{d^k}\)

				\item \label{tab:3:3}%
					\(\tilde s^{k}=\sum_{i=1}^N{\frac{1}{\gamma_i}\nabla h_i(u^k) - \frac{1}{N} \nabla f_i(u^k)}\)

				\item \label{tab:3:4}%
					\(y^k \in \Ts(\tilde s^k)\)%

				\item \label{tab:3:5}%
					\textbf{if}
					\(
						\sum_{i=1}^N  p_i(y^k,u^k)
						>
						\sum_{i=1}^N \tfrac{N}{\gamma_i}\DB_{h_i}(y^k,u^k)
					\)
					\textbf{then}

					\quad
					\(\gamma_i \gets \sigma \gamma_i\) for all \(i \in [N]\),
					update $\bar s^k$, and go to step \ref{tab:2:1}%

				\item
					\textbf{if}
					\(\mathcal{L}(y^k,u^k)\leq\mathcal{L}(v^k,z^k)\)

					\quad go to step \ref{tab:3:sts}%

				\item
					\textbf{else if}
					\(q_k = q_{\max}\)
					\textbf{then}

					\quad
					\(u^k = v^k\) and go to step \ref{tab:3:3}%

				\item
					\textbf{else}

					\quad
					\(\tau_k \gets \beta \tau_k\),~
					\(q_k \gets q_k + 1\),~
					and go to step \ref{tab:3:2}%
				\end{enumerate}

				\State \label{tab:3:sts}%
					\(s^k\gets\tilde s^k\)
			\end{algorithmic}

			\hrule

			\begin{algorithmic}[1]
			\addtocounter{ALG@line}{8}
			\State \label{tab:4}%
				\begin{enumerate}[label=\alph*:,ref=\ref*{tab:4}.\alph*,widest*=4]
				\item \label{tab:4:1}%
					\(\ttz[i^\ell]{k} \in \Ts(s^k)\)

				\item
					\textbf{if}
					\(p_{i^\ell}(\ttz[i^\ell]{k}, u^k) > \tfrac{N}{\gamma_{i^\ell}}\DB_{h_{i^\ell}}(\ttz[i^\ell]{k},u^k)\)
					\textbf{then}

					\quad
					\(\gamma_{i^\ell} \gets \sigma \gamma_{i^\ell}\),~
					update $s^k$,~
					and go to step \ref{tab:4:1}%
				\end{enumerate}
			\end{algorithmic}

			\hrule
			\hrule
			\label{tab:adaptive}
		\end{table}

		\subsection{Sparse Phase Retrieval with Squared Loss}\label{sec:phaseRetrieval}

			In the initial set of simulations, we evaluate the performance of SPIRAL on the sparse phase retrieval problem, which involves signal recovery based on intensity measurements. This problem finds applications in various fields, such as electron microscopy, speech recognition, optical imaging, and X-ray crystallography \cite{sun2018geometric,candes2015phase}. The sparse phase retrieval problem is formulated as follows:

			\begin{equation}
			\textstyle
			\minimize_{z \in \R^n} \frac{1}{N} \sum_{i=1}^{N} \textrm{Loss}(b_i,\langle a_i,z \rangle^2) + g(z),
			\label{eq:phase_retrieval}
			\end{equation}
			where $a_i \in \R^n$ and $b_i \in \R_+$ are data vectors. The objective of this optimization problem is to find a sparse vector $z \in \R^n$ that best approximates $b_i$ as $\langle a_i,z \rangle^2$ for all $i \in [N]$. Sparsity is introduced to account for noise and outliers, and this can be achieved by selecting nonsmooth functions such as the $\ell_0$ or $\ell_1$ norms. In this simulation, we consider the squared loss function $\textrm{Loss}(x,y)=\frac{1}{4}(x-y)^2$ and $g =\lambda \|\cdot\|_1$.
			To cast \eqref{eq:phase_retrieval} into the original optimization problem form \eqref{eq:problem_formulation}, we define the functions as follows:
			\[
				\textstyle
				f_i(z) = \frac{1}{4}(\langle a_i, z \rangle^2 - b_i)^2, ~\text{with}~
				h_i(z) = h(z) = \frac{1}{4} \|z\|^4 + \frac{1}{2} \|z\|^2.
			\]
			Although the cost functions $f_i$ do not possess Lipschitz continuous gradients, they exhibit smoothness relative to the reference function $h$ \cite[Lem. 5.1, Prop. 5.1 and 5.2]{bolte2018first}.

\begin{figure}[t]
				\centering
				\subcaptionbox{Performance of different algorithms.\label{fig:bregman:a}}
				{\includegraphics[width=0.4\linewidth]{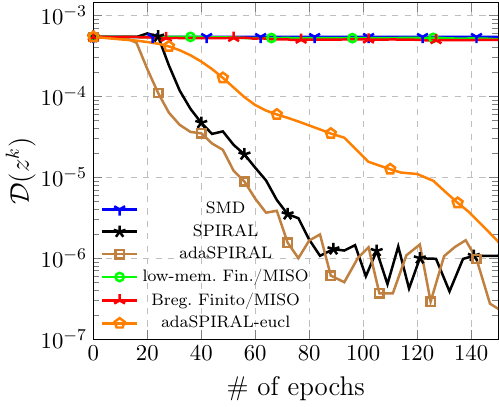}}
				\hspace{1cm}
				\subcaptionbox{Image recovery of adaSPIRAL.\label{fig:bregman:b}}
				{%
					\includegraphics[scale = 2]{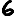}
					\hspace{0.5cm}
					\includegraphics[scale = 2]{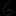}
					\hspace{0.5cm}
					\includegraphics[scale = 2]{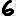}
					\vspace{0.5cm}
				}%
				\caption{Performance for the phase retrieval problem \eqref{eq:phase_retrieval} on a digit 6 image with $N=1280$, $n=256$. Image recovery is after $100$ epochs, including the original image (left), initialization (center), and output (right).}
				\label{fig:bregman}
			\end{figure}

			In this simulation, we consider $16 \times 16$ gray-scale images of digits from the dataset \cite{hastie2001elements}.
			The images are vectorized, so $n=256$. The matrix $A \in \R^{N \times n}$ with $a_i$ being its $i$th row, with $N=d \times n$ and $d=5$ is generated according to the procedure described in \cite{duchi2019solving}.
			We form this matrix as $A=[MS_1,\dots,MS_N]$ with $M \in \R^{n \times n}$ a normalized Hadamard matrix and $S_i$ diagonal sign matrices with the diagonal elements in $\{-1,1\}$. For noiseless data, $d=3$ is sufficient for a complete recovery. The measurements are corrupted by setting $b_i = 0$ with probability $p_c = 0.02$. Also, we set $\lambda=\frac{1}{N}$ by the hyperparameter search to have a visually good solution.
			Furthermore, the algorithms are initialized with the initialization scheme suggested in \cite{duchi2019solving}, and they converge to the same local optimal point.
			The performance of different algorithms is shown in \cref{fig:bregman:a} for a digit 6 image.

			As depicted in the figure, SPIRAL demonstrates significantly faster performance compared to the other algorithms. Even though the cost function in this scenario does not possess a Lipschitz continuous gradient, we evaluate the performance of adaSPIRAL, both in Bregman and Euclidean versions, in \cref{fig:bregman:a}.
			It is important to note that adaSPIRAL-eucl, implemented according to \cref{alg:LFinito:Euclidean} with the additional steps outlined in \cref{tab:adaptive} using dgfs $h_i = \frac{1}{2}\|\cdot\|^2$, does not require any prior knowledge of Lipschitz constants $L_i$.
			Remarkably, adaSPIRAL-eucl performs well on cost functions without Lipschitz continuous gradients, as verified by this simulation. Consequently, adaSPIRAL-eucl demonstrates potential applicability to a wider range of cost functions in various applications.
			Additionally, adaSPIRAL outperforms SPIRAL due to its ability to employ larger stepsizes that are dynamically updated as needed, thereby speeding up convergence.
			Furthermore, as shown in \cref{fig:bregman:b}, adaSPIRAL-eucl exhibits good image recovery capabilities even in highly corrupted initial conditions.

		\subsection{\texorpdfstring{$\ell_1$}{l1} Regularized Least Squares Problem}

			In this section, we evaluate the performance of SPIRAL for the Lasso problem, which is a convex optimization problem commonly used for regression tasks. The Lasso formulation is given by
			\begin{equation}
			    \textstyle
			    \minimize_{z \in \R^n} \frac{1}{2}\|Az-b\|^2_2 + \lambda\|z\|_1
			    \label{eq:lasso}
			\end{equation}
			where $A$ is a matrix with data vectors $a_i \in \R^n$ as its rows, and $b$ is a vector with corresponding labels $b_i \in \R$. The datasets used for regression tasks include the \emph{mg}, \emph{cadata}, \emph{housing}, and \emph{triazines} datasets obtained from LIBSVM. Additionally, synthetic datasets are generated using the procedure described in \cite[\S6]{nesterov2013gradient} for two different dimensions.
			The parameter $\lambda$ is appropriately set for each dataset.

			\begin{figure*}[t]
			    \centering
			    \includegraphics[width=\linewidth]{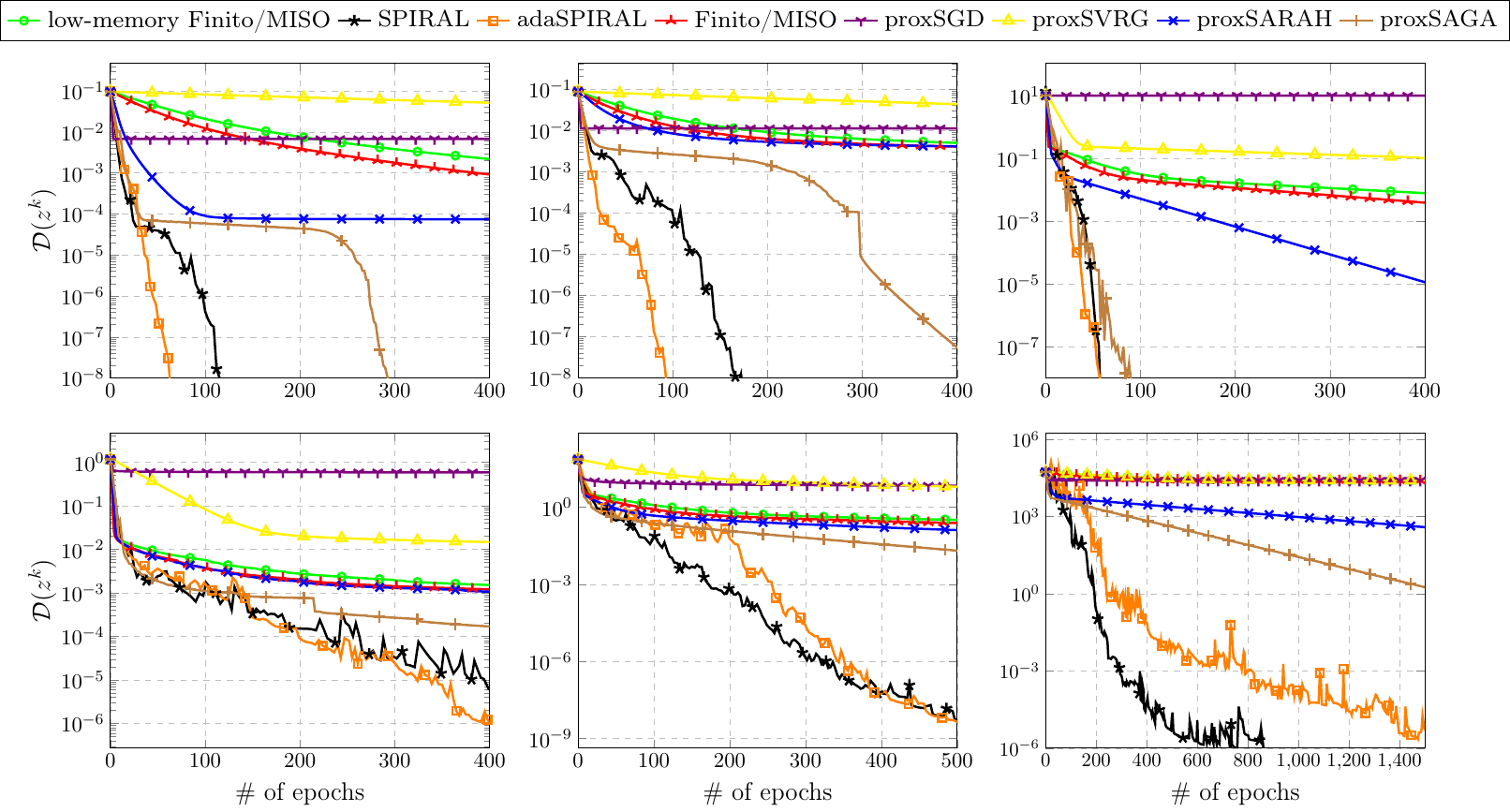}%
			    \caption[]{%
			        Performance of different algorithms for the Lasso problem \eqref{eq:lasso}. Synthetic dataset (top left) with $N=10000$, $n=400$, synthetic dataset (top center) with $N=300$, $n=600$, mg (top right) with $N=1385$, $n=6$, triazines (bottom left) with $N=186$, $n=60$, housing (bottom center) with $N=506$, $n=13$, and cadata (bottom right) with $N=20640$, $n=8$.}
			    \label{fig:lasso}%
			\end{figure*}

			\cref{fig:lasso} provides a comparison of different algorithms on six datasets. It is evident from the results that both SPIRAL and adaSPIRAL exhibit superior convergence performance compared to other algorithms, regardless of whether the datasets are synthetic or practical.
			Also the same speed up by adaSPIRAL is evident for most of the datasets.
			Note that adaSPIRAL does not require a priori knowledge of Lipschitz constants $L_i$, and still its performance is comparable with that of SPIRAL.
			Compared to Finito/MISO, low-memory Finito/MISO is worse, however, it does not need a large memory to store the gradient vectors.
			It is also observed that SPIRAL is particularly fast on dense datasets.

		\subsection{Nonnegative Principal Component Analysis}

			In this section, we investigate the problem of nonnegative principal component analysis (NN-PCA), which has also been studied in previous works \cite{reddi2016proximal,pham2020proxsarah}. The problem is formulated as follows:
			\begin{equation}\label{eq:nnpca}
				\minimize_{z \in \R^n}f(z)\coloneqq-\tfrac{1}{2N} \textstyle\sum_{i=1}^{N} z^T (a_ia_i^T) z
			\quad
				\stt \|z\| \leq 1,\ z \geq 0,
			\end{equation}
			where $N$ denotes the number of data points represented by $a_i \in \R^n$. To cast the problem with the form of \eqref{eq:problem_formulation}, we define $f_i(z)=-\frac{1}{2} z^T (a_ia_i^T) z$ and $g(z)=\delta_{\mathcal{B}}(z)$, where $\mathcal{B}\coloneqq \{w \in \R^n \mid \|w\| \leq 1, w \geq 0\}$ represents the constraints of the NN-PCA problem \eqref{eq:nnpca}.
			It is worth mentioning that SPIRAL allows for different smoothness constants and individualized stepsizes $\gamma_i$ for each of the functions $f_i$. In this case, the data points $a_i$ are not normalized to improve the output of NN-PCA analysis, providing a better representation of the dataset.
			All the algorithms employed for this nonconvex optimization problem are initialized with running proxSGD for 10 epochs, starting from the same initial point, ensuring convergence to a similar local optimum. The optimality criterion \eqref{eq:optimality_criteria} is reported as a function of the number of epochs.

			\begin{figure*}[t]
				\centering
				\includegraphics[width=\linewidth]{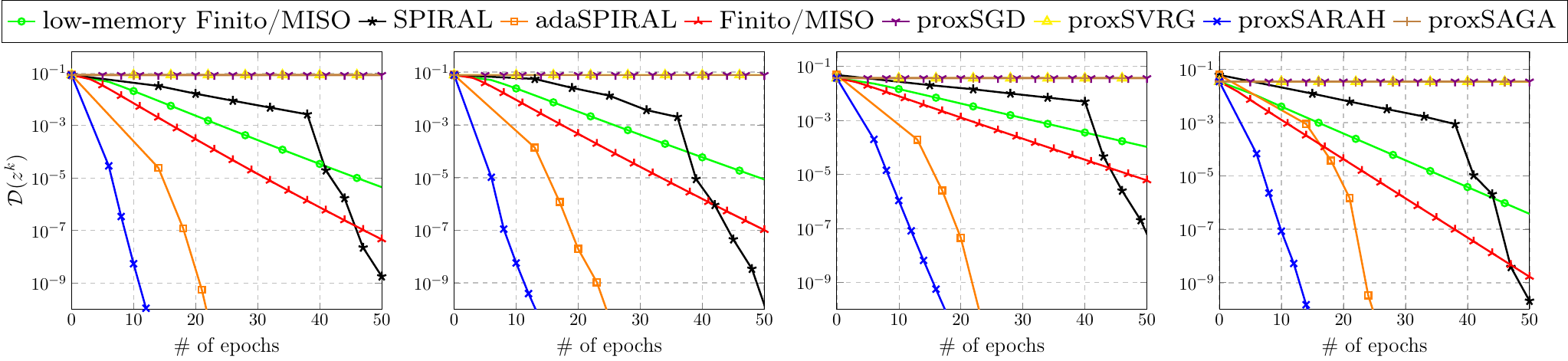}%
				\caption[]{%
					Performance of different algorithms for the NN-PCA problem of \eqref{eq:nnpca}. MNIST (left) with $N=60000$, $n=784$, covtype (left center) with $N=581012$, $n=54$, a9a (right center) with $N=32561$, $n=123$, and aloi (right) with $N=108000$, $n=128$.
				}
				\label{fig:nnpca}
			\end{figure*}

			As depicted in \cref{fig:nnpca}, the quasi-Newton updates in SPIRAL significantly enhance the convergence rate compared to (low-memory) Finito/MISO, which lacks such updates. Although proxSARAH exhibits faster convergence for this problem, it performs slower in the Lasso problem and is unable to handle non-Lipschitz differentiable cost functions, as demonstrated in the problem discussed in \cref{sec:phaseRetrieval}.

			In order to demonstrate how competitive different versions of SPIRAL are with state-of-the-art methods used in ML, also in terms of CPU time, performance comparisons are conducted in Appendix \ref{sec:cpu}.
			It is worth noting that despite SPIRAL's approximation of second-order information, the adopted quasi-Newton method, namely L-BFGS, demonstrates efficiency by relying solely on inexpensive level 1 BLAS operations, such as inner products, scalar multiplications, and additions.

	\section{Conclusion}\label{sec:con}

		This paper introduced SPIRAL, an optimization algorithm designed for solving regularized finite sum minimization problems.
		SPIRAL operates in a nonconvex setting and does not rely on the typical assumption of Lipschitz differentiability.
		Many existing methods that utilize quasi-Newton directions in finite sum settings either impose restrictive conditions or only achieve local convergence. In contrast, we demonstrated that SPIRAL achieves a superlinear convergence rate while ensuring \emph{global} convergence, all without the need for diminishing step sizes, and under standard mild assumptions.
		This is achieved by the introduction of a straightforward yet effective linesearch which is \emph{smart}, in the sense that it will never be triggered close enough to (sufficiently regular) solutions---an aspect validated also in our simulations.
		Moreover, it is observed that while addressing nonsmooth nonconvex problems, SPIRAL is still competitive with the state-of-the-art on classical convex problems, such as regularized least squares problems.
		Promising future research directions include the adaptation of SPIRAL to domains like distributed and federated learning.

	\begin{appendix}
		\section{Preliminaries}\label{sec:auxiliary:pre}

			\begin{fact}[basic properties \cite{chen1993convergence,nesterov2018introductory}]\label{lem:3point}%
				The following hold for a dgf \(H:\R^n\to\R\), \(x,y,z\in\R^n\):%
				\begin{factenumerate}[label=(\roman*),ref=\ref*{lem:3point}(\roman*),widest*=3]
				\item \label{thm:3P}%
					(three-point inequality)
					\(
						\DB_H(x,z)=\DB_H(x,y)+\DB_H(y,z)+\langle {x-y}, {\nabla H(y)-\nabla H(z)}\rangle.
					\)
					\cite[Lem. 3.1]{chen1993convergence}.
				\end{factenumerate}
				For any convex set \(\U\subseteq\R^n\) and \(u,v\in\U\) the following hold \cite[Thm. 2.1.5, 2.1.10]{nesterov2018introductory}:
				\begin{factenumerate}[resume,label=(\roman*),ref=\ref*{lem:3point}(\roman*),widest*=3]
				\item\label{thm:hstrconvex}%
					If \(H\) is \(\mu_{H,\U}\)-strongly convex on \(\U\), then
					\(
						\frac{\mu_{H,\U}}{2}\|v-u\|^2
					\leq
						\DB_{H}(v,u)
					\leq
						\frac{1}{2\mu_{H,\U}}\|\nabla H(v)-\nabla H(u)\|^2
					\).

				\item\label{thm:hC11}%
					If \(\nabla H\) is \(\ell_{H,\U}\)-Lipschitz on \(\U\), then
					\(
						\frac{1}{2\ell_{H,\U}}\|\nabla H(v)-\nabla H(u)\|^2
					\leq
						\DB_{H}(v,u)
					\leq
						\frac{\ell_{H,\U}}{2}\|v-u\|^2
					\).%
				\end{factenumerate}
			\end{fact}

			In the following, some properties of the Bregman Moreau envelope are highlighted. The interested reader is referred to \cite{ahookhosh2021bregman} and \cite{kan2012moreau} for proofs and further properties.

			\begin{fact}[\protect{Basic properties of $\FBEG$ and $\prox_\FBEG$, \cite{ahookhosh2021bregman,kan2012moreau}}]\label{prop:Moreau_properties}%
				Let $H:\R^n\to \R$ denote a dgf (cf. \cref{def:legendre_kernel}), and $\phi:\R^n\to \Rinf$ be a proper lsc, and lower bounded function. Then, the following hold:
				\begin{factenumerate}[label=(\roman*),ref=\ref*{prop:Moreau_properties}(\roman*),widest*=3]
				\item \label{prop:prox_properties}%
					$\prox_\FBEG$ is locally bounded, compact-valued, and outer semicontinuous;

				\item \label{prop:Moreau_properties_2}%
					$\FBEG$ is finite-valued and continuous; it is locally Lipschitz if so is $\nabla H$;

				\item \label{prop:Moreau_properties_1}%
					$\FBEG(z) = \phi(v) + \DB_{H}(v, z) \leq \phi(y) + \DB_{H}(y, z)$ with any $y, z\in \R^n$, $v \in \prox_\FBEG(z)$. Hence, $\FBEG(z) \leq \phi(z)$;

				\item \label{prop:Moreau_properties_3}%
					$\inf \phi = \inf \FBEG$ and $\argmin \FBEG =  \argmin \phi$;

				\item \label{prop:Moreau_properties_4}%
					$\FBEG$ is level-bounded iff so is $\phi$.
				\end{factenumerate}
			\end{fact}

			The following fact studies sufficient conditions for Lipschitz continuity of the Bregman proximal mapping and continuity of the Moreau envelope, both of which are crucial to the theory developed in \cref{lem:descent_lemma,thm:superlinear}.

			\begin{fact}[\protect{\cite[Lem. A.2]{latafat2021bregman}}]\label{lem:reg:env}%
				Let \(\mathcal{V}_i\subseteq\R^n\) be nonempty and convex, \(i\in[N]\), and let \(\mathcal{V}\coloneqq\mathcal{V}_1\times\cdots\times\mathcal{V}_N\).
				Additionally to \cref{ass:basic}, suppose that \(g\) is convex, and \(h_i\), $i\in[N]$, is \(\ell_{h_i}\)-smooth and \(\mu_{h_i}\)-strongly convex on \(\mathcal{V}_i\).
				Then, the following hold for function \(\hat{H}\) as in \eqref{eq:H_hat_def} with \(\gamma_i\in(0,\nicefrac{N}{L_{f_i}})\), \(i\in[N]\):
				\begin{factenumerate}[label=(\roman*),ref=\ref*{lem:reg:env}(\roman*),widest*=2]
				\item \label{lem:TLip}%
					\(\prox_\FBE\) is $\bar L$-Lipschitz continuous on \(\mathcal{V}\) for some constant $\bar L\geq0$.
				\end{factenumerate}
				If in addition \(f_i\) and \(h_i\) are twice continuously differentiable on \(\mathcal{V}_i\), \(i\in[N]\), then
				\begin{factenumerate}[resume,label=(\roman*),ref=\ref*{lem:3point}(\roman*),widest*=2]
				\item \label{lem:C1+}%
					\(\FBE\) is continuously differentiable on \(\mathcal{V}\) with \(\nabla\FBE=\nabla^2\hat{H}\circ(\id-\prox_\FBE)\).
				\end{factenumerate}
			\end{fact}

			The following fact establishes the equivalence between problems \eqref{eq:problem_formulation} and \eqref{eq:recasted_problem_formulation}.

			\begin{fact}[\protect{\cite[Lem. A.1]{latafat2021bregman}}]\label{lem:problem_equivalence}%
				Let the functions $\varphi$ and $\Phi$ be as in \eqref{eq:problem_formulation} and \eqref{eq:recasted_problem_formulation}, respectively.
				Then,
				\begin{factenumerate}[label=(\roman*),ref=\ref*{lem:problem_equivalence}(\roman*),widest*=2]
				\item \label{lem:equiv:1}%
					$\partial \Phi(\bm x) = \{\bm v = (v,\dots,v) \mid \sum_iv_i \in \partial \varphi(x)\}$ if $\bm x=(x,\dots,x) \in \Delta$, and is empty otherwise.%

				\item \label{lem:equiv:2}%
					$\Phi$ has the KL property at $\bm x = (x,\dots,x)$ iff so does $\varphi$ at $x$.
					In this case, the desingularizing functions are the same up to a positive scaling.
				\end{factenumerate}
			\end{fact}

		\section{Omitted lemmas}\label{sec:omitted_lemmas}

			\begin{lemma}\label{lem:Tri}
				Suppose that \cref{ass:basic,ass:superlinear} hold and that $\varphi$ is level bounded. Consider the sequence generated by \cref{alg:LBFinito}.
				Then, for every $\ell\in[N]$ there exists $c_\ell >0$ such that
				\begin{equation}\label{eq:xu}
					\|\ttz[\ell]{k} - u^k\| \leq c_\ell \|\ttz[1]{k} - u^k\|.
				\end{equation}
				\begin{proof}
					By level boundedness of $\varphi$ and \cref{lem:descent_lemma}, $(\bm u^k)_{k\in\N}$, $(\bm z^k)_{k\in\N}$, $(\ttz[i]{k})_{k\in\N}, i\in[N]$ are contained in a nonempty bounded set $\bm \U$.
					By \cref{ass:fh}, $h_i$ is locally strongly convex and locally Lipschitz, which along with \cref{ass:g:cvx} and \cref{lem:reg:env} implies that $\prox_\FBE$ is $\bar L$-Lipschitz on a convex subset of $\bm \U$ for some $\bar L > 0$.
					Without loss of generality and for the sake of simplicity, we assume the cyclic sweeping rule in the incremental loop, i.e., $i^\ell=\ell$.
					Note that the following proof can be easily cast into the case of cyclic sweeping without replacement.
					Arguing by induction, for $\ell=1$, \eqref{eq:xu} holds trivially.
					Suppose that the claim holds for some $\ell\geq 1$. Then, by triangular inequality and the definition of $\tbz[\ell]{k}$ in step \ref{algLift:1:t3} of \cref{tab:lifted}
					\begin{align*}
						\|\ttz[\ell+1]{k} - u^k\|
					={} &
						\tfrac{1}{\sqrt{N}}\|\tbz[\ell+1]{k} - \bm u^k\|
					\\
					\leq{} &
						\tfrac{1}{\sqrt{N}} \|\tbz[1]{k} - \bm u^k\| + \tfrac{1}{\sqrt{N}} \|\tbz[\ell+1]{k} - \tbz[1]{k}\|
					\\
						\text{\footnotesize Lip. continuity of \(\prox_\FBE\) and \cref{tab:lifted}}
					\leq{} &
						\tfrac{1}{\sqrt{N}} \|\tbz[1]{k} - \bm u^k\| + \tfrac{\bar L}{\sqrt{N}}\|\tz[\ell+1]{k} - \bm u^k\|
					\\
						\text{\footnotesize \eqref{eq:tz_def}}
					\leq{} &
						\textstyle\|\ttz[1]{k} - u^k\| + \tfrac{\bar L}{\sqrt{N}}\sum_{j\leq \ell}\|\ttz[j]{k} - u^k\|
					\\
						\text{\footnotesize (induction)}
					\leq{} &
						\underbracket[0.5pt]{\textstyle \left(1 + \tfrac{\bar L}{\sqrt{N}} \sum_{j\leq \ell} c_j \right)}_{\coloneqq c_{\ell+1}} \|\ttz[1]{k} - u^k\|,
					\end{align*}
					establishing \eqref{eq:xu}.
				\end{proof}
			\end{lemma}

			\begin{lemma}\label{lem:Tri:2}
				In addition to the assumptions in \cref{lem:Tri}, suppose that the directions $d^k$ in step \ref{alg:1:dir} satisfy \(\|d^k\|\leq D\|z^k-v^k\|\)
				for some \(D\geq 0\).
				Then, $\|z^{k+1}- z^k\| \leq C\|\bm z^k - \tz[N]{k-1}\|$ holds for some positive $C$.
				\begin{proof}
					By the same reasoning as in \cref{lem:Tri}, \(\prox_\FBE\) is $\bar L$-Lipschitz continuous on a bounded convex set containing the iterates $(\bm u^k)_{k\in\N}$, $(\bm z^k)_{k\in\N}$, $(\ttz[i]{k})_{k\in\N}, i\in[N]$.
					It follows from the assumption on $\|d^k\|$ and step \ref{algLift:2:ls:1} of \cref{tab:lifted} that
					\begin{equation}\label{eq:zu}
						\|\bm z^k - \bm u^k\| \leq (1-\tau_k)\|\bm z^k - \bm v^k\|  + \tau_k \|\bm d^k\| \leq (1-\tau_k + \tau_k D) \|\bm z^k - \bm v^k\| \leq \eta_1 \|\tz[N]{k-1} - \bm z^k\|,
					\end{equation}
					where $\eta_1 = \bar L(1-\tau_k + \tau_k D)$ and Lipschitz continuity of the proximal mapping was used in the last inequality.
					Further using triangular inequality yields
					\begin{align*}
					\numberthis\label{eq:12}
						\|\bm u^k - \tz[N]{k-1}\|
					\leq{} &
						\|\bm u^k - \bm z^k\| + \|\bm z^k - \tz[N]{k-1}\|
					\\
						\text{\footnotesize \eqref{eq:zu}}
					\leq{}&
						\left(\eta_1 + 1\right)\|\bm z^k - \tz[N]{k-1}\|, \quad \text{and}
					\\
						\|\ttz[1]{k} - u^k\|
					={}&
						\tfrac1{\sqrt{N}}\|\tbz[1]{k} - \bm u^k\|
					\\
					\leq{}&
						\tfrac1{\sqrt{N}}\|\tbz[1]{k} - \bm z^k\| + \tfrac1{\sqrt{N}} \|\bm z^k - \bm u^k\|
					\\
						\text{\footnotesize Lip. continuity of \(\prox_\FBE\) and \cref{tab:lifted}}
					\leq{}&
						\tfrac{\bar L}{\sqrt{N}}\|\bm u^k - \tz[N]{k-1}\| + \tfrac1{\sqrt{N}} \|\bm z^k - \bm u^k\|
					\\
						\text{\footnotesize \eqref{eq:12}, \eqref{eq:zu}}
					\leq{}&
						\tfrac1{\sqrt{N}}\big((\bar L+1)\eta_1 + \bar L\big)\|\bm z^k - \tz[N]{k-1}\|.
					\end{align*}
					Using this along with triangular inequality yields
					\[
						\textstyle \|\tz[N]{k} - \bm u^k\| = \sum_{\ell=1}^N \|\ttz[\ell]{k} - u^k\|
					\stackrel{\text{\eqref{eq:xu}}}{\leq}
						\sum_{\ell=1}^N  c_\ell \|\ttz[1]{k} - u^k\| \leq \eta_2\|\bm z^k - \tz[N]{k-1}\|,
					\]
					where $\eta_2 = \sum_{\ell=1}^N \tfrac{c_\ell}{\sqrt{N}}\big((\bar L+1)\eta_1 + \bar L\big)$.
					This inequality combined with \eqref{eq:12} yields
					\begin{align*}
						\|z^{k+1}- z^{k}\|
					={} &
						\frac{1}{\sqrt{N}}\|\prox_\FBE(\tz[N]{k}) - \prox_\FBE(\tz[N]{k-1})\|
					\leq
						\frac{\bar L}{\sqrt{N}} \|\tz[N]{k}- \tz[N]{k-1}\|
					\\
					\leq{} &
						\frac{\bar L}{\sqrt{N}} \|\tz[N]{k}- \bm u^k\|+ \frac{\bar L}{\sqrt{N}} \|\bm u^k - \tz[N]{k-1}\|
					\leq
						\frac{\bar L}{\sqrt{N}} (\eta_1 + \eta_2 + 1)\|\bm z^k - \tz[N]{k-1}\|.
					\end{align*}
					The claimed inequality follows from Lipschitz continuity of $\prox_\FBE$ and
					the inclusion in step \ref{algLift:1:t1} of \cref{tab:lifted}.
				\end{proof}
			\end{lemma}

		\section{Omitted proofs}\label{sec:auxiliary:glb}

			\phantomsection
			\addcontentsline{toc}{subsection}{Proof of Theorem \ref*{thm:global}}%
			\subsection*{Proof of \cref{thm:global}}\label{aux:thm:global}
				By level boundedness of $\varphi$ and \cref{lem:descent_lemma}, $(\bm u^k)_{k\in\N}$, $(\bm z^k)_{k\in\N}$, $(\ttz[i]{k})_{k\in\N}$ are contained in a nonempty convex bounded set $\bm \U$, where owing to \cref{ass:fh}, $h_i$ and consequently $\hat{H}$ are strongly convex.
				It then follows from \cref{thm:hstrconvex}, \cref{lem:descent_lemma:res}, and \cref{lem:Tri:2} that
				\(
					\|z^{k+1} - z^k\| \to 0.
				\)
				Therefore, the set of limit points of \((z^k)_{k\in\N}\) is nonempty compact and connected \cite[Rem. 5]{Bolte2014Proximal}.
				By \cref{lem:descent_lemma:4,lem:descent_lemma_stationary} the limit points are stationary for \(\varphi\), and $\FBE(\bm z^k) = \mathcal{L}(v^k, z^k) \to \varphi_\star$.
				In the trivial case \(\FBE(\bm z^k) = \mathcal{L}(v^k, z^k) = \varphi_\star\) for some $k$, the claims follow from \cref{lem:descent_lemma}.
				Assume that $\FBE(\bm z^k)> \varphi_\star$ for $k\in \N$. The KL property for $\Phi$ is implied by that of $\varphi$ due to \cref{lem:problem_equivalence}, with desingularizing function $\psi(s)=\rho s^{1-\theta}$ with exponent $\theta \in (0,1)$.
				Let $\Omega$ denote the set of limit points of $(\bm z^k=(z^k,\ldots, z^k))_{k\in\N}$. Since  $\hat{H}$ is  strongly convex, \cite[Lem. 5.1]{yu2021kurdyka} can be invoked to infer that the function $\MM(\bm w, \bm x) = \Phi(\bm w) + \DB_{\hat{H}}(\bm w, \bm x)$ also has the KL property with exponent $\nu \in \max\{\theta,\frac{1}{2}\}$ at every point  \((\bm z^\star, \bm z^\star)\) in the compact set \(\Omega \times \Omega\). Moreover, by \eqref{eq:H_hat_def} \(\MM(\bm z^\star,\bm z^\star) = \Phi(\bm z^\star)= \varphi_\star\) where \cref{lem:descent_lemma:4} was used in the last equality.
				Recall that $\bm z^{k} \in \prox_\FBE(\tz[N]{k-1})$ as in step \ref{algLift:1:t1} of \cref{tab:lifted}.
				Therefore, $ \partial \MM(\bm z^{k}, \tz[N]{k-1}) = (\underbrace{\partial \Phi(\bm z^k) + \nabla \hat{H}(\bm z^k) - \nabla \hat{H}(\tz[N]{k-1})}_{ \ni 0, \; \text{by  \eqref{thm:OC} in the lifted space}}, \nabla^2 \hat{H}(\tz[N]{k-1})(\tz[N]{k-1} - \bm z^{k}) )$, resulting in    %
				\begin{equation}
					\dist(0,\partial \MM(\bm z^{k}, \tz[N]{k-1})) \leq \|\nabla^2 \hat{H}(\tz[N]{k-1})\|\|\tz[N]{k-1} - \bm z^k\| \leq c \|\tz[N]{k-1} - \bm z^k\|
					\label{eq:dist_part_M}
				\end{equation}
				where $c = \sup_k\|\nabla^2 \hat{H}(\tz[N]{k-1})\| > 0$ is finite due to \(\ttz[N]{k}\) being bounded (cf. \cref{lem:descent_lemma:bounded}) and continuity of $\nabla ^2 \hat{H}$.
				Considering \eqref{eq:KL_def} with \eqref{eq:dist_part_M}, since \(\MM(\bm z^k, \tz[N]{k-1}) = \FBE(\tz[N]{k-1}) \to \varphi\) from above, and that \((\bm z^k, \tz[N]{k-1})_{k\in\N}\) is bounded and accumulates on $\Omega\times \Omega$, up to discarding iterates the following holds
				\begin{equation}\label{eq:psi_dev}
					\textstyle
					\psi'\big(\FBE(\tz[N]{k-1})- \varphi_\star\big)
				=
					\psi'\big(\MM(\bm z^{k}, \tz[N]{k-1})- \MM(\bm z^\star, \bm z^\star)\big)  \geq \frac{1}{c\|\tz[N]{k-1} - \bm z^k\|},
				\end{equation}
				where $\psi=\rho s^{1-\nu}$ is a desingularizing function for $\MM$ on $\Omega \times \Omega$.
				Let us define
				\begin{equation}
					\Delta_k \coloneqq \psi(\FBE(\tz[N]{k-1})-\varphi_\star)
				=
					\rho [\FBE(\tz[N]{k-1})-\varphi_\star]^{1-\nu}
				\leq
					\rho [\rho (1-\nu)c\|\tz[N]{k-1} - \bm z^k\|]^{\frac{1-\nu}{\nu}}.
					\label{eq:delta}
				\end{equation}
				Then,
				\(
					\Delta_k^{\frac{\nu}{1-\nu}} \leq c \rho^{\frac{1}{1-\nu}}(1-\nu)\|\tz[N]{k-1} - \bm z^k\|.
				\)
				Concavity of $\psi$ also implies
				\begin{equation}
					\Delta_k - \Delta_{k+1}
				\geq
					\psi'(\FBE(\tz[N]{k-1})-\varphi_\star)(\FBE(\tz[N]{k-1})-\FBE(\tz[N]{k}))
				\stackrel{\text{\ref{eq:psi_dev}}}{\geq}
					\frac{\FBE(\tz[N]{k-1}) - \FBE(\tz[N]{k})}{c \|\tz[N]{k-1} - \bm z^k\|}.
					\label{tmp:superlinear_1}
				\end{equation}
				On the other hand by \eqref{eq:des_lem_2} and \eqref{eq:des_lem_1}
				\begin{equation}
					\FBE(\tz[N]{k+1}) - \FBE(\tz[N]{k}) \leq - \DB_{\hat{H}}(\bm z^{k+1}, \tz[N]{k})
					\leq -\frac{\mu_{\hat{H}}}{2}\|\bm z^{k+1} - \tz[N]{k}\|^2,
					\label{tmp:superlinear_2}
				\end{equation}
				where \cref{thm:hstrconvex} was used and $\mu_{\hat{H}}$ denotes its strong convexity modulus.
				Combining \eqref{tmp:superlinear_1} and \eqref{tmp:superlinear_2},
				\begin{equation}
					\textstyle
					\Delta_k - \Delta_{k+1}
				\geq
					\eta\|\tz[N]{k-1} - \bm z^k\|
				\geq
					\frac{\eta}{C} \|\bm z^{k+1} - \bm z^k\|,
					\label{eq:delta_diff}
				\end{equation}
				with some constant $\eta > 0$ where the last inequality follows from \cref{lem:Tri:2}. Hence, $(\|\bm z^{k+1} - \bm z^k\|)_{k\in\N}$ has finite length and is thus convergent.
				It then follows from \cref{lem:descent_lemma_stationary} that $(\bm z^k)_{k\in\N}$ converges to a stationary point of $\varphi$.
				Combining \eqref{eq:delta} and \eqref{eq:delta_diff} we have,
				\begin{equation}
					\Delta_{k+1} \leq \Delta_k - \alpha \Delta_k^{{\frac{\nu}{1-\nu}}}
					\label{eq:delta_rel}
				\end{equation}
				with some appropriate $\alpha \geq 0$. Hence, if $\nu=\frac{1}{2}$, i.e. $\theta \in (0,\frac{1}{2}]$ for $\Phi$, in \eqref{eq:delta_rel} we have $\Delta_{k+1} \leq (1 - \alpha) \Delta_k$. As $\alpha > 0$ and $\frac{\Delta_{k+1}}{\Delta_k} > 0$, then $(1-\alpha) \in (0,1)$ concluding $\Delta_k$ is Q-linearly convergent to zero.
				By \eqref{eq:delta} we then conclude $(\FBE(\tz[N]{k-1}))_{k\in\N}$ is convergent Q-linearly and by \cref{prop:Moreau_properties_1}, where we have $\varphi(z^k)=\Phi(\bm z^k) \leq \FBE(\tz[N]{k-1})$, we conclude $(\varphi(z^k))_{k\in\N}$ is convergent R-linearly.
				Moreover, the inequality \eqref{eq:delta_diff} implies that $(\|z^{k+1} - z^k\|)_{k\in\N}$ is R-linearly convergent, thus so is $(z^k)_{k\in\N}$.

		\section{CPU time}\label{sec:cpu}
\begin{figure}
				\centering
				\includegraphics[width=0.38\linewidth]{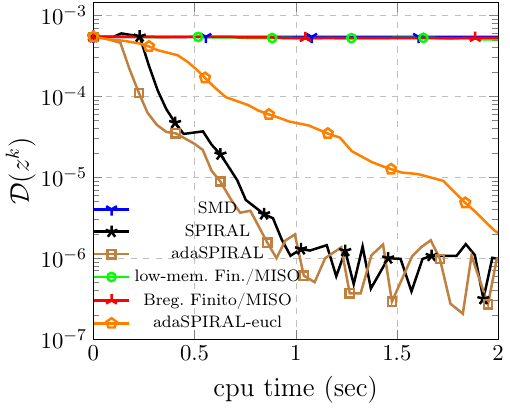}%
				\captionsetup{font={footnotesize}, labelfont={footnotesize}}
			    \caption[]{%
			        Performance of different algorithms versus cpu time on the phase retrieval problem \eqref{eq:phase_retrieval} for $550$ epochs on a digit 6 image with $N=1280$, $n=256$.%
			    }%
			    \label{fig:bregman:cpu}%
			\end{figure}

\begin{figure}
				\centering
				\includegraphics[width=\linewidth]{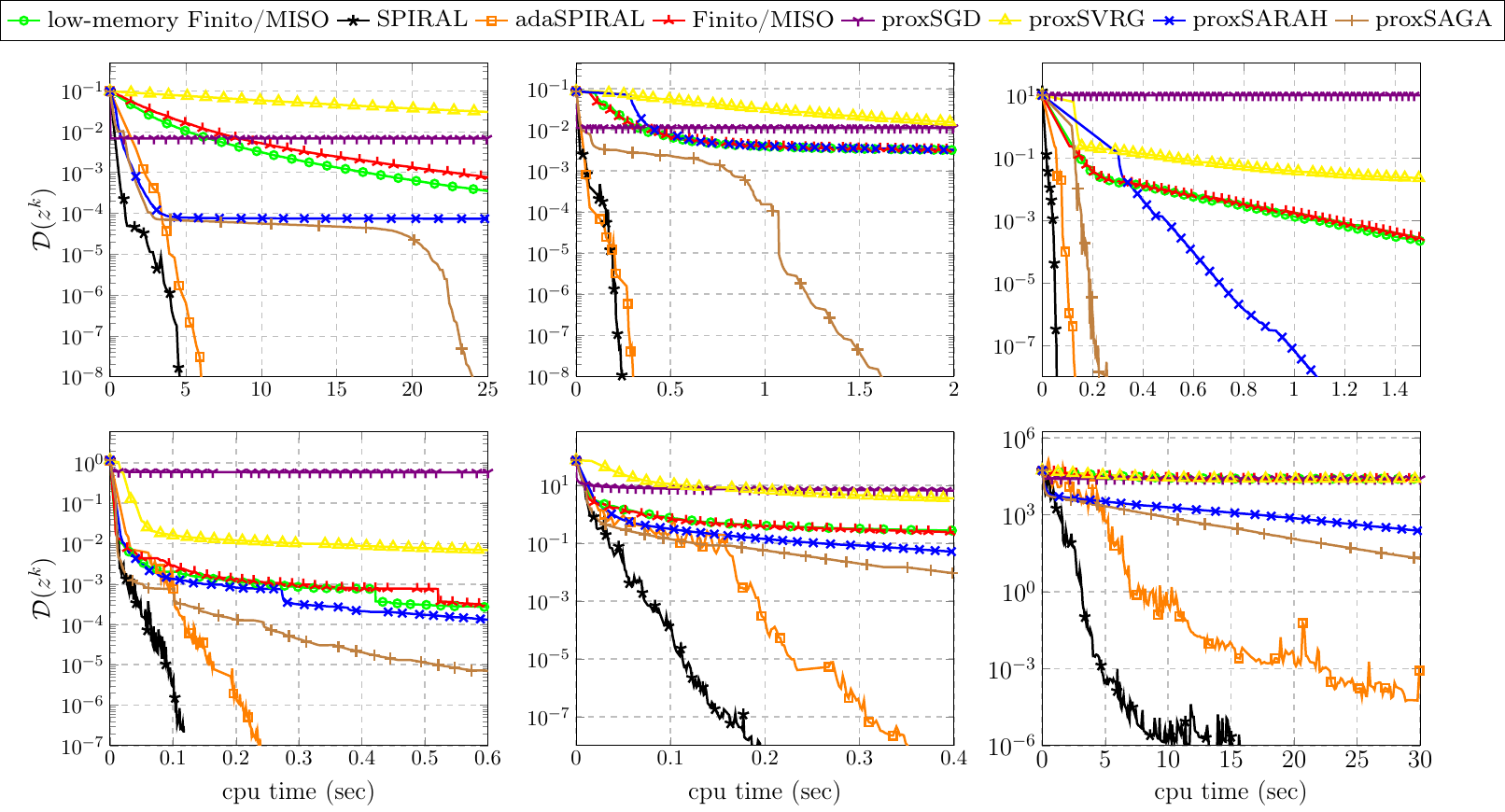}%
				\captionsetup{font={footnotesize}, labelfont={footnotesize}}
			    \caption[]{%
			    	Performance of different algorithms versus cpu time on the lasso problem of \eqref{eq:lasso} for $50$ epochs. Synthetic dataset (top left) with $N=10000$, $n=400$, synthetic dataset (top center) with $N=300$, $n=600$, mg (top right) with $N=1385$, $n=6$, triazines (bottom left) with $N=186$, $n=60$, housing (bottom center) with $N=506$, $n=13$, and cadata (bottom right) with $N=20640$, $n=8$.
			    }%
			    \label{fig:lasso:cpu}%
			\end{figure}

\begin{figure}
			    \centering
				\includegraphics[width=\linewidth]{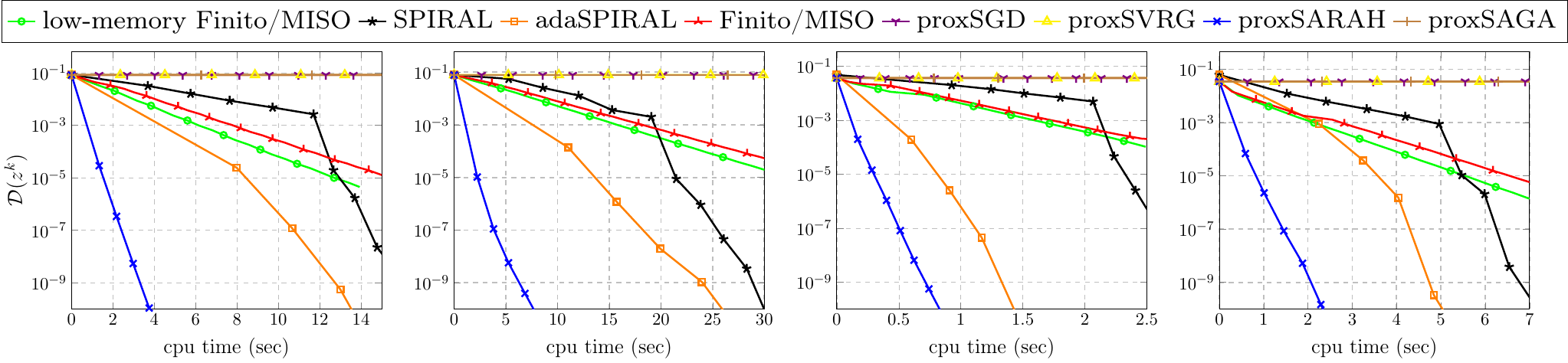}%
				\captionsetup{font={footnotesize}, labelfont={footnotesize}}
			    \caption[]{%
			        Performance of different algorithms versus cpu time on the NN-PCA problem of \eqref{eq:nnpca} for $500$ epochs. MNIST (left) with $N=60000$, $n=784$, covtype (left center) with $N=581012$, $n=54$, a9a (right center) with $N=32561$, $n=123$, and aloi (right) with $N=108000$, $n=128$.
			    }
			    \label{fig:nnpca:cpu}
			\end{figure}

			The performance results presented in \cref{sec:sim} are also reported versus CPU time.
			According to the numerical comparisons in \cref{fig:lasso:cpu,fig:bregman:cpu,fig:nnpca:cpu}, the proposed algorithm features relatively cheap iterations and has comparable computational complexity per epoch compared to the other algorithms.


		\section{Algorithm variants}

			\subsection{Adaptive variant}\label{sec:auxiliary:adaptive}

				In this section, the implementation of \cref{tab:adaptive} is further discussed.
				In \cref{tab:adaptive}, for the first iterate, i.e. $k=0$, the vectors $\ttz[i]{-1}$ are initially considered equal to $z^{\rm init}$ for all $i\in[N]$.
				Also, note that the linesearch in step \ref{tab:3:5} of \cref{tab:adaptive} backtracks to step \ref{tab:2:1}, rather than step \ref{tab:3:4}.
				Performing the linesearches in this intertwined fashion is observed to result in acceptance of good directions and reduction in the overall computational complexity \cite{demarchi2021proximal,pas2021alpaqa}.
				We refer the reader to \cite{demarchi2021proximal} for the theoretical justification for the effectiveness of this procedure.
				Note that in \cref{alg:LFinito:Euclidean}, in the Euclidean case, the same backtrackings can be used with dgfs $h_i = \frac{1}{2} \|\cdot\|^2$.
				The backtracking linesearches in the first block of \cref{tab:adaptive} do not require storing $\ttz[i]{k}$ and can be performed efficiently.
				In step \ref{tab:1:2} $\textstyle \sum_{i=1}^N p_i(\cdot,\ttz[i]{k})$ may be evaluated by storing the scalars $\sum_{i=1}^N f_i(\ttz[i]{k})$ and $\sum_{i=1}^N \langle \nabla f_i(\ttz[i]{k}), \ttz[i]{k} \rangle$ and one vector $\sum_{i=1}^N \nabla f_i(\ttz[i]{k}) \in \R^n$  while performing step \ref{alg:forloop:end} of the algorithm. Similar tricks apply to the computation of the Bregman distances, functions $p_i$ in other backtracking linesearches of \cref{tab:adaptive}, and updating the vectors $s^k, \bar s^k$, and $\tilde s^k$.

			\subsection{Euclidean variant}\label{sec:Euclidean}

			\begin{algorithm}[t]
				\caption{SPIRAL - Euclidean version}
				\label{alg:LFinito:Euclidean}

				\begin{algorithmic}[1]
				\item[\textbf{Require}]
					\begin{tabular}[t]{@{}l@{}}
						\(z^{\rm init}\in\R^n\),~
						$\Gamma=\{\gamma_1,\dots,\gamma_N\}$ with \(\gamma_i\in(0,\nicefrac{N}{L_{i}})\), \(i\in[N]\),~
						\(\beta\in(0,1)\),
					\\
						\(\hat\gamma\coloneqq\bigl(\sum_{i=1}^N\gamma_i^{-1}\bigr)^{-1}\),~
						\(s^{0}=z^{\rm init}-\frac{\hat \gamma}{N}\sum_{i=1}^N\nabla f_i(z^{\rm init})\),
					\\
						maximum number of backtracks \(q_{\max} \in \N\cup\{\infty\}\) {\footnotesize(e.g. \(q_{\max}=2\))},~
						\(K\in\N\)
					\end{tabular}

				\item[\textbf{Repeat for} $k=0,1,\dots,K$]
					\State
						\(z^k\in\prox_{\hat\gamma g}(s^k)\)

					\State
						$\bar s^k = z^k-\tfrac{\hat{\gamma}}{N}\sum_{i=1}^{N}\nabla f_{i}(z^k)$
					\Comment{(full update)}

					\State
						\(v^k\in\prox_{\hat\gamma g}(\bar s^k)\)

					\State \label{alg:2:dir}%
						choose $d^k\in\R^n$ at $z^k$ {\footnotesize(e.g. based on a quasi-Newton method for solving $\res(z)=0$)}

					\State \label{alg:2:ls}%
						set $\tau_k = 1, q_k=0$
					\Comment{(linesearch)}
						\begin{enumerate}[label=\alph*:,ref=\ref*{alg:2:ls}.\alph*]
						\item \label{alg:2:ls:1}
							\(u^k=\tau_k z^k+(1-\tau_k)v^k+\tau_k{d^k}\)

						\item \label{alg:2:ls:2}%
							\(\tilde s^k = u^k-\tfrac{\hat{\gamma}}{N}\sum_{i=1}^{N}\nabla f_{i}(u^k)\)
						\Comment{(full update)}

						\item \label{alg:2:ls:3}%
							\(y^k\in\prox_{\hat \gamma g}(\tilde s^k)\)

						\item \label{eq:2:ls_condition}%
							\textbf{if}
							\(\mathcal{L}(y^k,u^k)\leq\mathcal{L}(v^k,z^k)\)

							\quad
							go to step \ref{alg:2:sts}%

						\item
							\textbf{else if}
							\(q_k = q_{\rm max}\)
							\textbf{then}

							\quad
							\(u^k = v^k\),~
							\(\tilde s^k = u^k-\tfrac{\hat{\gamma}}{N}\sum_{i=1}^{N}\nabla f_{i}(u^k)\),~
							and go to step \ref{alg:2:sts}

						\item
							\textbf{else}

							\quad
							\(\tau_k \gets \beta \tau_k\),~
							\(q_k \gets q_k + 1\),~
							and go to step \ref{alg:2:ls:1}%
						\end{enumerate}

					\State \label{alg:2:sts}%
						$s^k \gets \tilde s^k$

					\For{ \(\ell = 1,\ldots, N\) }
					\Comment{(incremental loop)}
						\State
							randomly choose \(i^\ell\in[N]\) without replacement

						\State
							\(\tilde z^k_{i^\ell}\in\prox_{\hat\gamma g}(s^k)\)

						\State
							\(
								s^k
							\gets
								s^k+\tfrac{\hat{\gamma}}{N}\left(\nabla f_{{i^\ell}}(u^k)-\nabla f_{{i^\ell}}(\tilde z^k_{i^\ell})\right)
								+
								\tfrac{\hat{\gamma}}{\gamma_{{i^\ell}}}\left(\tilde z^k_{i^\ell} - u^k\right)
							\)
					\EndFor

					\State
						$s^{k+1} \gets s^k$

				\item[\textbf{Return} $z^K$]
				\end{algorithmic}
			\end{algorithm}

			In this section, the proposed algorithm in the Euclidean version is outlined in \cref{alg:LFinito:Euclidean}, when the functions $f_i$ have Lipschitz continuous gradients with constants $L_i$.
			In this case, the distance generating functions are $h_i=\frac{1}{2}\|\cdot\|^2$ and consequently the Bregman distances are simplified to $\DB_{h_i}(y,x)=\frac{1}{2}\|y-x\|^2$.
	\end{appendix}


	\clearpage
	\phantomsection
	\addcontentsline{toc}{section}{References}%
	\bibliographystyle{plain}%
	\bibliography{Bibliography.bib}


\end{document}